\documentclass[preprint,12pt]{elsarticle}

\usepackage{graphicx}
\usepackage{multirow}
\usepackage{amsmath,amssymb,amsfonts}
\usepackage{amsthm}
\usepackage{mathrsfs}
\usepackage{mathtools}
\usepackage[title]{appendix}
\usepackage{xcolor}
\usepackage{textcomp}
\usepackage{manyfoot}
\usepackage{booktabs}
\usepackage{algorithm}
\usepackage{algorithmicx}
\usepackage{algpseudocode}
\usepackage{listings}
\usepackage{geometry}

\theoremstyle{definition}
\newtheorem{defn}{Definition}[section]

\newtheorem{rem}[defn]{Remark}
\theoremstyle{plain}
\newtheorem{lem}[defn]{Lemma}
\newtheorem{thm}[defn]{Theorem}
\newtheorem{cor}[defn]{Corollary}
\newtheorem{pr}[defn]{Proposition}
\newtheorem{Assumption}{Assumption}[section]

\def\Ec{{\cal E}}
\def\Fc{{\cal F}}

\def\Lc{{\cal L}}
\def\Mc{{\cal M}}

\def\Pc{{\cal P}}

\renewcommand{\P}{\mathbb{P}}
\newcommand{\e}{\mathbb{E}}
\newcommand{\R}{\mathbb{R}}
\newcommand{\rom}[1]{
  \textup{\uppercase\expandafter{\romannumeral#1}}
}
\def \ito{Itô's formula }

\theoremstyle{thmstyleone}

\theoremstyle{thmstyletwo}

\theoremstyle{thmstylethree}

\journal{arXiv}

\begin{document}

\begin{frontmatter}

\title{Multi-dimensional reflected McKean-Vlasov BSDEs with the obstacle depending on both the first unknown and its distribution}

\author[1]{Ruisen Qian\corref{cor1}}
\ead{rsqian19@fudan.edu.cn}
\cortext[cor1]{Corresponding author}
\affiliation[1]{ organization={School of Mathematical Sciences, Fudan University}, 
                 city={Shanghai},
                 postcode={200433}, 
                 country={China}
 }

\begin{abstract}
    The paper studies a multi-dimensional mean-field reflected backward stochastic differential equation (MF-RBSDE) with a reflection constraint depending on both the value process $Y$ and its distribution $[Y]$. 
    We establish the existence, uniqueness and the stability of the solution of MF-RBSDE. 
    We also investigate the associated interacting particle systems of RBSDEs and prove a propagation of chaos result. Lastly, we investigate the relationship between MF-RBSDE and an obstacle problem for partial differential equations in Wasserstein space within a Markovian framework. 
    Our work provides a connection between the work of Briand et al. (2020) on BSDEs with normal reflection in law and the work of Gegout-Petit and Pardoux (1996) on classical multi-dimensional RBSDEs. 
\end{abstract}


\end{frontmatter}

\section{Introduction}

Since the seminal paper by Pardoux and Peng \cite{PP90}, backward stochastic differential equations (BSDEs) have been extensively studied. 
The solution of a typical BSDE with a random driver $f$ and a terminal condition $\xi$ is a pair of progressively measurable processes $(Y, Z)$ that satisfy the dynamic 
\begin{equation}\label{BSDE}
	Y_t=\xi+\int_t^T f(s,Y_s,Z_s)ds -\int_t^T Z_sdB_s, \quad 0\leq t\leq T.
\end{equation}

Pardoux and Peng \cite{PP90} demonstrated the existence and uniqueness of such a solution, and since then, many extensions of the dynamic have been proposed and studied.

El Karoui et al. \cite{EKPPM97} introduced the notion of a reflected BSDE. The solution of a reflected BSDE contains an additional adapted non-decreasing process $K$ with $K_0=0$, such that the triplet of processes $(Y,Z,K)$ satisfies the dynamic
\begin{equation}\label{RBSDE}
	\begin{aligned}
		Y_t=\xi +\int_t^T f(s,Y_s,Z_s)ds-\int_t^T Z_s dB_s + K_T-K_t, \quad 0\leq t\leq T
	\end{aligned}
\end{equation}
with a chosen constraint on the solution. El Karoui et al. \cite{EKPPM97} considered the constraint of the form
\begin{equation}
	\begin{aligned}
		Y_t\geq S_t, \quad 0\leq t\leq T,
	\end{aligned}
\end{equation}
with a continuous progressively measurable obstacle process $S$, 
along with a Skorokhod condition 
\begin{equation}
	\begin{aligned}
		\int_0^T (Y_t-S_t) dK_t=0
	\end{aligned}
\end{equation}
to ensure the minimality of the solution.
The multi-dimensional case where the process $Y$ is constrained in a convex domain is studied by Gegout-Petit and Pardoux \cite{GPP96}.

More recently, Briand et al. \cite{BCDH20,BEH18,BH21} proposed and studied a class of BSDEs with normal reflection in law, where the distribution $\mu$ of the $Y$ component of the solution is required to satisfy the constraint and the corresponding Skorokhod condition 
\begin{equation}\label{El}
	\begin{aligned}
		h(\mu_t)\geq0, \quad 0\leq t\leq T, \,\,\mathrm{and}\,\, \int_0^T h(\mu_s) dK_s=0, 
	\end{aligned}
\end{equation}
for some Lions differentiable concave functional $h:\Pc_2(\R^n)\rightarrow \R$. The dynamic for this class of BSDEs writes
\begin{equation}\label{RBSDE-BCDH20}
	\begin{aligned}
		Y_t=\xi +\int_t^T f(s,Y_s,Z_s)ds-\int_t^T Z_sdB_s + \int_t^T \partial_\mu h(\mu_s)(Y_s) dK_s, \quad 0\leq t\leq T. 
	\end{aligned}
\end{equation}
Briand et al. \cite{BCDH20} showed that the system (\ref{El})-(\ref{RBSDE-BCDH20}) admits a unique solution if $K$ is only allowed to be deterministic. This type of reflected BSDEs finds applications in quantile hedging problems \cite{FL99,FL00} and super-hedging problems under risk measure constraints \cite[Section 6]{BEH18}. 

The aim of this paper consists in enlarging the result of \cite{BCDH20}, and providing a connection to the classical results of reflected BSDEs \cite{GPP96}. 
We study a general class of mean-field reflected BSDEs where the constraint takes the form 
\begin{equation}\label{refl}
	H(Y_t,\mu_t)\geq 0, 
\end{equation}
where the functional $H$ depends on both the value process $Y$ and its distribution $\mu$. Moreover, the driver $f$ is allowed to depend on the variables $(Y, Z)$ as well as their joint distribution $\nu$. 
The dynamic of our problem writes 
\begin{equation}
	Y_t=\xi+\int_t^T f(s,Y_s,Z_s,\nu_s)ds - \int_t^T Z_s dB_s + \int_t^T \partial_y H(Y_s, \mu_s)dK_s +\widetilde{\e}\left[\int_t^T \partial_\mu H(\widetilde{Y}_s,\mu_s)(Y_s)d\widetilde{K}_s\right] ,
\end{equation}
for all $t\in [0,T]$. 
For example, let $H(y,\mu)=y+h(\mu)$, then the reflection part of this dynamic becomes 
\begin{equation}
	K_T - K_t + \int_t^T \partial_\mu h(\mu_s)(Y_s) d\e[K_s] ,
\end{equation}
which can be viewed as a combination of the reflections in \eqref{RBSDE} and \eqref{RBSDE-BCDH20}. 

Reflected BSDEs with constraints depending on both the value process $Y$ and its distribution $\mu$ have found applications in insurance and risk management. One such application is in the pricing of guaranteed life endowment policies with a surrender/withdrawal option, as described by \cite{DEH19}. In this context, the obstacle $h$ contains a bonus option, which is linked to the distribution of the possible surplus realized by the average of all involved contracts. 

The paper establishes the well-posedness of mean-field reflected BSDEs, connects these problems with the mean-field limit of a system of reflected BSDEs, and investigates the relationship between the mean-field reflected BSDEs and an obstacle problem for partial differential equations (PDEs) in Wasserstein space. 

The rest of the paper is organized as follows. In Section \ref{sec-framework}, the problem is formulated in detail and the assumptions are clarified. We also give some a priori estimates of the solution, which will be used later on.

In Sections \ref{sec-unique-stability} and \ref{sec-exist}, the well-posedness of the solution is analyzed. We first present a stability result. Next, we establish the uniqueness of the solution by combining the stability result and an analysis of the reflection component. Finally, we prove the existence of a solution using the penalization technique.

In Section \ref{sec-limit}, we attempt to interpret the mean-field reflected BSDEs at the particle level. We consider a corresponding particle system and investigate the limiting properties of the solution of this system. We demonstrate that the solution of the mean-field reflected BSDEs is the mean-field limit of the particle system.

Finally, in Section \ref{sec:obstacle problem}, we establish a connection between our mean-field reflected BSDEs and an obstacle problem for PDEs in Wasserstein space. We show that, given that the problem is formulated within a Markovian framework, the solution of mean-field reflected BSDEs provides a probabilistic representation of a viscosity solution of an obstacle problem in Wasserstein space.

\smallskip

\paragraph*{Notations.}

Throughout this paper, we will work on a classical Wiener space $(\Omega, \cal F, \mathbb{P})$ with a finite time horizon $T>0$ and a $d$-dimensional standard Brownian motion $B=(B_t)_{0\leq t \leq T}$. We endow the probability space $(\Omega, \cal F, \mathbb{P})$ with the filtration $\mathbf{F}=(\Fc_t)_{0\leq t\leq T}$ generated by the Brownian motion $B$. We also denote by:
\begin{itemize}
	\item $S^{2,n}$ the set of $R^n$-valued continuous adapted processes $Y$ on $[0, T]$ such that
	$\\ \left \| Y \right \|_{S^{2,n}}^2 := \e \left[ \sup_{t\in [0, T]} |Y_t|^2\right]<\infty$,
	\item $H^{2,n}$ the set of $R^{n\times d}$-valued predictable processes $Z$ such that $\left \| Z \right \|_{H^{2,n}}^2 := \e \left[\int_0^T |Z_t|^2 dt\right]<\infty$,
	\item $A^{2,1}$ the subset of $S^{2,1}$ consisting of non-decreasing processes starting from $0$,
	\item $\Pc_2(\R^n)$ the Wasserstein space equipped with 2-Wasserstein distance $W_2(\cdot,\cdot)$, 
	\item $[\xi]$ the distribution of a random variable $\xi$. 
\end{itemize}

\section{Framework}\label{sec-framework}

\subsection{Mean-field reflected BSDEs}

Throughout this paper, we consider the following problem: for all $t\in[0,T]$,
\begin{equation}\label{MF-RBSDE}
	Y_t=\xi+\int_t^T f(s,Y_s,Z_s,\nu_s)ds - \int_t^T Z_s dB_s + \int_t^T \partial_y H(Y_s, \mu_s)dK_s +\widetilde{\e}\left[\int_t^T \partial_\mu H(\widetilde{Y}_s,\mu_s)(Y_s)d\widetilde{K}_s\right] ,
\end{equation}
\begin{equation}\label{refl}
	H(Y_t,\mu_t)\geq 0,\quad \int_0^T H(Y_s, \mu_s)dK_s=0 ,
\end{equation}
where $\mu_s:=[Y_s], \nu_s:=[(Y_s,Z_s)]$ denote the distribution of $Y_s$ and the joint distribution of $(Y_s,Z_s)$, and $\widetilde{Y}_s, \widetilde{K}_s$ are independent copies of $Y_s,K_s$. 

We define a solution of (\ref{MF-RBSDE})-(\ref{refl}) as a triple of progressively measurable processes $(Y,Z,K)$ with values in $\R^n\times\R^{n\times d}\times\R$ such that $K$ is continuous, non-decreasing, and $K_0=0$. 

Some special cases of our problem are: 
\begin{enumerate}
	\item $H(y,\mu)=H(y)$, meaning the constraint depends solely on the value process $Y$. In this case, the system (\ref{MF-RBSDE})-(\ref{refl}) corresponds to the classical reflected BSDE with normal reflection, as studied in \cite{GPP96}.
	\item $H(y,\mu)=H(\mu)$, meaning the constraint depends only on the distribution $\mu$ of the value process $Y$. In this case, the system (\ref{MF-RBSDE})-(\ref{refl}) corresponds to the reflected BSDE with normal constraint in law, which is investigated in \cite{BCDH20,BEH18}. We note that the reflection terms of (\ref{MF-RBSDE}) reduce to $\int_t^T \partial_\mu H(\mu_s)(Y_s)\,\e\left[ dK_s \right]$. Consequently, the reflection process $K$ is usually assumed to be deterministic in order to obtain uniqueness results. 
	\item $H(y,\mu)=y+h(\mu)$ for some concave functional $h: \Pc_2(\R)\to\R$, meaning the constraint decomposes into a part that depends on $y$ and a part that depends on $\mu$. This is a generalization of cases 1 and 2, which is also described earlier in the Introduction. 
\end{enumerate}

We study the system (\ref{MF-RBSDE})-(\ref{refl}) under the following assumptions:
\begin{Assumption}\label{H}
	The driver $f$, the functional $H$, and the terminal condition $\xi$ satisfy:
	\begin{enumerate}
		\item\label{H_f} 
		$f$ is a mapping from $\Omega \times [0, T] \times \R^n \times \R^{n\times d} \times \Pc_2(\R^n\times\R^{n\times d})$ into $\R^n$ such that 
		\begin{enumerate}
			\item\label{H_f_H2} The process $(f(t,0,0,\delta_0))_{0\leq t\leq T}$ is progressively measurable and
			\begin{equation*}
				\e\left[\int_0^T |f(t,0,0,\delta_0)|^4 dt \right]<\infty. 
			\end{equation*}
			\item There exists a constant $L \geq 0$, such that for all $t\in[0,T], i\in\{1,2\}$, and for all $y_i \in \R^n$, $z_i\in \R^{n*d}$, $\nu_i\in\Pc_2(\R^n\times\R^{n\times d})$,
			\begin{equation*}
				|f(t, y_1,z_1,\nu_1)-f(t, y_2,z_2,\nu_2)| \leq L\left(|y_1-y_2|+|z_1-z_2|+W_2(\nu_1,\nu_2)\right).
			\end{equation*}
		\end{enumerate}
		\item\label{H_H}
		The functional $H: \R^n \times \Pc_2(\R^n) \to \R$ has the following properties in $\R^n$ 
		\begin{enumerate}
			\item\label{H-differentiable-y}
			For any $\mu\in \Pc_2(\R^n)$, the function $\R^n\ni y\mapsto H(y,\mu)$ is twice continuously differentiable in $\R^n$, and the functions $\partial_y H$ and $\partial^2_{yy} H$ are jointly continuous in $(y,\mu)$.

			\item\label{H-differentiable-mu}
			For any $y\in\R^n$, the functional $\Pc_2(\R^n)\ni \mu\mapsto H(y,\mu)$ is continuously L-differentiable and, for any $\mu\in\Pc_2(\R^n)$, there exists a version of the function $\R^n\ni v\mapsto\partial_\mu H(y,\mu)(v)$, such that the functional $\partial_\mu H$ is jointly continuous in $(y,\mu,v)$.

			\item\label{H-differentiable-mu-v}
			For the version $\partial_\mu H$ mentioned above and for any $(y,\mu,v)\in\R^n\times\Pc_2(\R^n)\times\R^n$, the function $\R^n\ni v\mapsto\partial_\mu H(y,\mu)(v)$ is continuously differentiable in $\R^n$, and its derivative, denoted by $\R^n\ni v\mapsto \partial_v\partial_\mu H(y,\mu)(v)\in\R^{n\times n}$, is jointly continuous in $(y,\mu,v)$.

			\item\label{H-Lipschitz}
			There exist $0<\beta<M<\infty$, such that for all $(y,\mu,v)\in R^n\times\Pc_2(\R^n)\times\R^n$,
			\begin{equation}\label{H-beta}
				|\partial_y H(y,\mu)|\geq \beta,
			\end{equation}
			and
			\begin{equation}\label{H-M}
				|\partial_y H(y,\mu)|+|\partial^2_{yy} H(y,\mu)|+|\partial_\mu H(y,\mu)(v)|+|\partial_v\partial_\mu H(y,\mu)(v)|\leq M.
			\end{equation}

			\item\label{Hy-Lipschitz}
			The functionals $\partial_y H(y,\cdot)$ and $\partial_\mu H(y,\cdot)(\cdot)$ are Lipschitz continuous, i.e., 
			for all $y\in\R^n$, $Y_1,Y_2\in L^2(\Omega, \R^n)$, and $\mu_1,\mu_2\in\Pc^2(\R^n)$, 
			\begin{equation}
			\begin{gathered}\label{Hy-L}
				\left| \partial_y H(y,\mu_1)- \partial_y H(y,\mu_1) \right|
				\leq 
				L \, W_2(\mu_1,\mu_2), \\
				\e\left| \partial_\mu H(y,[Y_1])(Y_1)- \partial_\mu H(y,[Y_2])(Y_2) \right|^2
				\leq 
				L \, \e\left| Y_1-Y_2 \right|^2. 
			\end{gathered}
			\end{equation}

			\item\label{Hc} For all $y_1,y_2,v\in\R^n$ and for all $\mu_1,\mu_2\in\Pc_2(\R^n)$,
			\begin{equation}
			\begin{gathered}\label{Hx-Hy}
				\partial_y H(v,\mu_1)\cdot \partial_\mu H(y_2,\mu_2)(v)\geq 0, \\
				\partial_\mu H(y_1,\mu_1)(v) \cdot\partial_\mu H(y_2,\mu_2)(v)\geq 0.
			\end{gathered}
			\end{equation}

			\item\label{H-concave}
			The functional $H$ is concave:
			for all $(y_1,\mu_1),(y_2,\mu_2)\in \R^n\times\Pc_2(\R^n)$,
			\begin{equation}\label{H1-H2}
				H(y_2,\mu_2)-H(y_1,\mu_1)
				\leq
				\partial_y H(y_1,\mu_1)\cdot (y_2-y_1)
				+\e \left[\partial_\mu H(y_1,\mu_1)(Y_1)\cdot(Y_2-Y_1)\right],
			\end{equation}
			whenever $Y_1$ and $Y_2$ are square integrable random variables with distributions $\mu_1$ and $\mu_2$. 

		\end{enumerate}
		\item
		The terminal value $\xi$ is an $\Fc_T$-measurable random variable with $\e|\xi|^4<\infty$, such that
		\begin{equation*}\label{H-xi-geq-0}
			H(\xi,\left[\xi\right])\geq 0. 
		\end{equation*}
	\end{enumerate}
\end{Assumption}

\begin{rem}
	\begin{enumerate}
		\item Assumption \ref{H} 2.(a)-(c) are equivalent to the Assumption(Joint Chain Rule) in \cite[Section 5.6.4]{CD18}, which is used to derive chain rules in the Wasserstein space. 
		\item We assume that $\int_0^T |f(t,0,0,\delta_0)|^4 dt$ 
		and $\xi$ have finite 4th moments. These technical assumptions are necessary in the proof of Lemma \ref{lem-Ym-Zm-Cauchy}, which is crucial for establishing the existence of a solution. 
	\end{enumerate}
\end{rem}

Let us begin with a simple proposition, which provides a pivotal point for obtaining a priori estimates. 
\begin{pr}\label{pr-widetilde-y}
	Under Assumption \ref{H}, there exists $\widehat{y}\in\R^n$, such that
	\begin{equation}
		\begin{aligned}
			H(\widehat{y},\delta_{\widehat{y}})\geq0.
		\end{aligned}
	\end{equation}
\end{pr}

\begin{proof}
	We define $(y(t))_{t\geq0}$ as the solution of
	\begin{equation}
		\begin{aligned}
			y(t)=\int_0^t \partial_y H(y(s),\delta_{y(s)})ds,
		\end{aligned}
	\end{equation}
	and derive from \eqref{H-beta} and \eqref{Hx-Hy} that
	\begin{equation}
		\begin{aligned}
			& H(y(t),\delta_{y(t)})-H(0,\delta_0)
			\\
			=    &
			\int_0^t \left( \left|\partial_y H(y(s),\delta_{y(s)})\right|^2 + \partial_\mu H(y(s),\delta_{y_(s)})(y(s)) \cdot \partial_y H(y(s),\delta_{y(s)})
			\right) ds \\
			\geq &
			\int_0^t  \left|\partial_y H(y(s),\delta_{y(s)})\right|^2
			ds         \\
			\geq &
			\beta^2 t.
		\end{aligned}
	\end{equation}
	Therefore, there exists $t^*\geq0$, such that $H(y(t^*),\delta_{y(t^*)})\geq0$. We set $\widehat{y}:=y(t^*)$ to obtain the desired result.
\end{proof}

\subsection{A priori estimate}

In the following, we provide some useful a priori estimates for any solution of the MF-RBSDE (\ref{MF-RBSDE})-(\ref{refl}). 

\begin{lem}\label{lem-a-priori-supE}
	Under Assumption \ref{H}. Let $(Y,Z,K)$ be a solution of (\ref{MF-RBSDE})-(\ref{refl}) in $S^{2,n}\times H^{2,n}\times A^{2,1}$. Then, there exists $C>0$, such that $(Y,Z)$ satisfies the following:
	\begin{equation}\label{est-supE-Y}
		\begin{aligned}
			\sup_{0\leq t \leq T}\e\left[ 	\left| Y_t \right|^2\right] +  \e\left[\int_0^T \left|Z_s\right|^2ds\right]\leq C \e \left[ |\xi|^2 + \int_0^T |f^0(s)|^2ds + |\widehat{y}|^2 \right],
		\end{aligned}
	\end{equation}
	\begin{equation}\label{est-supE-Y4}
		\begin{aligned}
			\sup_{0\leq t \leq T}\e\left[ 	\left| Y_t \right|^4\right] \leq C \e \left[ |\xi|^4 + \int_0^T |f^0(s)|^4ds + |\widehat{y}|^4 \right],
		\end{aligned}
	\end{equation}
	where $\widehat{y}$ is the constant vector appearing in Proposition \ref{pr-widetilde-y}, and $f^0(s):=f(s,0,0,\delta_0)$. 
\end{lem}

\begin{proof}

	Let $\Delta Y:=Y-\widehat{y}$ (recalling that $\widehat{y}\in\R^n$ and $H(\widehat{y},\delta_{\widehat{y}})\geq0$). Applying \ito to $e^{\alpha t}|\Delta Y_t|^2$, we obtain 
	\begin{equation}\label{ito10}
		\begin{aligned}
			  & e^{\alpha t} |\Delta Y_t|^2 + \int_t^T e^{\alpha s} |Z_s|^2ds                                                                   \\
			= &
			e^{\alpha T} |\xi-\widehat{y}|^2
			- \int_t^T \alpha e^{\alpha s} |\Delta Y_s|^2 ds
			+  \int_t^T  2 e^{\alpha s} \Delta Y_s \cdot f(s,Y_s,Z_s,\nu_s) ds                                                                        \\
			  & + \int_t^T 2 e^{\alpha s} \Delta Y_s \cdot \partial_y H(Y_s,\mu_s) dK_s
			+ \int_t^T 2 e^{\alpha s} \Delta Y_s \cdot \widetilde{\e}\left[ \partial_\mu H(\widetilde{Y}_s,\mu_s)(Y_s)d\widetilde{K}_s \right] \\
			  & - \int_t^T 2 e^{\alpha s} \Delta Y_s \cdot Z_s dB_s.                                                                             \\
		\end{aligned}
	\end{equation}

	From the Lipschitz continuity of $f$ and Young's inequality, we obtain 
	\begin{equation}
	\begin{aligned}
		& 2 e^{\alpha s} \Delta Y_s \cdot f(s,Y_s,Z_s,\nu_s) \\
		\leq & 
		2  e^{\alpha s} |\Delta Y_s| \left[ |f^0(s)| + L\left(|Y_s| + |Z_s| + W_2(\nu_s,\delta_{(0,0)})\right) \right]
		\\
		\leq & 
		2  e^{\alpha s} |\Delta Y_s| \left[ |f^0(s)| + L\left(|\Delta Y_s| + |\widehat{y}| + |Z_s| + W_2(\nu_s,\delta_{(\widehat{y},0)}) + |\widehat{y}| \right) \right]
		\\
		\leq & 
		e^{\alpha s} \left[  \left(1 +2L +12 L^2\right) |\Delta Y_s|^2 + \e|\Delta Y_s|^2 + |f^0(s)|^2 + |\widehat{y}|^2  +\frac{1}{4} |Z_s|^2 + \frac{1}{4} \e |Z_s|^2\right]. 
	\end{aligned}
	\end{equation}
	Setting $\alpha:= 2 +2L +12 L^2$ and taking the expectation on both sides of (\ref{ito10}), we obtain 
	\begin{equation}\label{eito10}
		\begin{aligned}
			     & \e \left[ e^{\alpha t} |\Delta Y_t|^2 +  \frac{1}{2} \int_t^T e^{\alpha s} |Z_s|^2ds \right]                                                        \\
			\leq &
			\e \left[ 	e^{\alpha T} 	|\xi-\widehat{y}|^2
			+  \int_t^T  e^{\alpha s} \left(\left| f^0(s) \right|^2 +|\widehat{y}|^2\right)ds \right]                                                                                             \\
			&+2\e \left[\int_t^T  e^{\alpha s} \Delta Y_s \cdot \left(
				\partial_y H(Y_s,\mu_s) dK_s
				+ \widetilde{\e}\left[ \partial_\mu H(\widetilde{Y}_s,\mu_s)(Y_s)d\widetilde{K}_s \right]
				 \right)\right]. \\
		\end{aligned}
	\end{equation}

	From the concavity of $H$, the inequality $H(\widehat{y},\delta_{\widehat{y}})\geq 0$, and Fubini's theorem, we obtain 
	\begin{equation}\label{est-K0}
		\begin{aligned}
			& \e \left[\int_t^T  e^{\alpha s} \Delta Y_s \cdot \left(
				\partial_y H(Y_s,\mu_s) dK_s
				+ \widetilde{\e}\left[ \partial_\mu H(\widetilde{Y}_s,\mu_s)(Y_s)d\widetilde{K}_s \right]
				 \right)\right]    \\
			=    & \e \left[\int_t^T  e^{\alpha s} \left(\Delta Y_s \cdot \partial_y H(Y_s,\mu_s) +  \widetilde{\e}\left[\Delta \widetilde{Y}_s \cdot \partial_\mu H(Y_s,\mu_s)(\widetilde{Y}_s)  \right] \right) dK_s \right]   \\
			& + \e \left[\int_t^T  e^{\alpha s} \left(   \Delta Y_s \cdot \widetilde{\e}\left[ \partial_\mu H(\widetilde{Y}_s,\mu_s)(Y_s)d\widetilde{K}_s \right]
			-\widetilde{\e}\left[\Delta \widetilde{Y}_s \cdot \partial_\mu H(Y_s,\mu_s)(\widetilde{Y}_s)  \right]  dK_s \right) \right] \\
			=    & \e \left[\int_t^T  e^{\alpha s} \left(\Delta Y_s \cdot \partial_y H(Y_s,\mu_s) +  \widetilde{\e}\left[\Delta \widetilde{Y}_s \cdot \partial_\mu H(Y_s,\mu_s)(\widetilde{Y}_s)  \right] \right) dK_s \right]   \\
			\leq & \e \left[\int_t^T  e^{\alpha s} \left(H(Y_s,\mu_s)-H(\widehat{y}, \delta_{\widehat{y}})\right) dK_s\right]        
			\leq 0,
		\end{aligned}
	\end{equation}
	with the last inequality coming from the Skorokhod condition and the inequality $ H(\widehat{y},\delta_{\widehat{y}}) \geq0$.

	From (\ref{eito10}) and (\ref{est-K0}), we have 
	\begin{equation}\label{est-supE0}
		\begin{aligned}
			\sup_{0\leq t \leq T} \e \left[ e^{\alpha t} |\Delta Y_t|^2\right]
			+ \e \left[  \int_0^T e^{\alpha t} |Z_s|^2ds \right]
			\leq & C \e \left[ e^{\alpha T} |\xi-\widehat{y}|^2
			+  \int_0^T  e^{\alpha s} \left| f^0(s)\right|^2 ds \right], \\
		\end{aligned}
	\end{equation}
	and prove \eqref{est-supE-Y}. Then, \eqref{est-supE-Y4} is deduced by 
    applying \ito to $e^{\alpha t}|\Delta Y_t|^4$, 
    \begin{equation}\label{eq:ito-4}
    \begin{aligned}
        e^{\alpha t}|\Delta Y_t|^4 +  \int_t^T 4 e^{\alpha s} |\Delta Y_s \cdot Z_s|^2 ds 
    = 
        e^{\alpha T} |\xi-\widehat{y}|^4
        - \int_t^T \alpha e^{\alpha s}|\Delta Y_s|^4 ds 
        - \int_t^T 2 e^{\alpha s}|\Delta Y_s|^2 d|\Delta Y_s|^2, 
    \end{aligned}
    \end{equation}
    and following a similar approach as previously described. 
\end{proof}

\begin{lem}\label{lem-a-priori-Esup}
	Let Assumption \ref{H} be satisfied. Let $(Y,Z,K)$ be a solution of (\ref{MF-RBSDE})-(\ref{refl}) in $S^{2,n}\times H^{2,n}\times A^{2,1}$. Then, there exists $C>0$, such that $(Y,Z)$ satisfies the following:
	\begin{equation}
		\begin{aligned}
			\e\left[\sup_{0\leq t \leq T}	\left| Y_t \right|^2 +  \int_0^T \left|Z_s\right|^2ds \right]
			\leq C \e\left[ |\xi|^2 + \int_0^T |f^0(s)|^2ds + |\widehat{y}|^2 + (K_T)^2 \right].
		\end{aligned}
	\end{equation}
\end{lem}

\begin{proof}
	Using the same notation as in Lemma \ref{lem-a-priori-supE}, we take the supremum in $t$ and the expectation in (\ref{ito10}) to obtain
	\begin{equation}\label{e16}
		\begin{aligned}
			& \e\left[\sup _{0 \leq t \leq T} e^{\alpha t}\left|\Delta Y_t\right|^2+\frac{1}{2} \int_0^T e^{\alpha s}\left|\Delta Z_s\right|^2 d s\right]
			\\
			\leq & \e\left[e^{\alpha T}|\xi-\widehat{y}|^2+\int_0^T e^{\alpha s}\left(|f^0(s)|^2+|\widehat{y}|^2\right) d s\right]+2 \e \sup _{0 \leq t \leq T}\left|\int_t^T e^{\alpha s} \Delta Y_s \cdot Z_s d B_s\right| \\
			& +\e \sup _{0 \leq t \leq T}\left[\int_t^T e^{\alpha s} \Delta Y_s \cdot\left(\partial_y H(Y_s, \mu_s) d K_s+\widetilde{\e}\left[\partial_\mu H(\widetilde{Y}_s, \mu_s)\left(Y_s\right) d \widetilde{K}_s\right]\right)\right].
			\\
		\end{aligned}
	\end{equation}

	From the concavity of $H$, the boundedness of $\partial_y H$ and $\partial_\mu H$, the equation $\e\left[ \int_0^T H(Y_t,\mu_t)dK_t\right]=0$, and \eqref{est-supE-Y}, we obtain  
	\begin{equation}\label{e16-K-L}
		\begin{aligned}
			& \e \sup _{0 \leq t \leq T}\left[\int_t^T e^{\alpha s} \Delta Y_s \cdot\left(\partial_y H(Y_s, \mu_s) d K_s+\widetilde{\e}\left[\partial_\mu H(\widetilde{Y}_s, \mu_s)\left(Y_s\right) d \widetilde{K}_s\right]\right)\right]
			\\
			\leq & \e \sup _{0 \leq t \leq T}\left[\int_t^T e^{\alpha s}\left(\Delta Y_s \cdot \partial_y H(Y_s, \mu_s)+\widetilde{\e}\left[\Delta \widetilde{Y}_s \cdot \partial_\mu H(Y_s, \mu_s)(\widetilde{Y}_s)\right]\right) d K_s\right] \\
			& +\e \sup _{0 \leq t \leq T}\left[-\int_t^T e^{\alpha s} \widetilde{\e}\left[\Delta \widetilde{Y}_s \cdot \partial_\mu H(Y_s, \mu_s)(\widetilde{Y}_s)\right] d K_s\right] \\
			& +\e \sup _{0 \leq t \leq T}\left[\int_t^T e^{\alpha s} \Delta Y_s \cdot \widetilde{\e}\left[\partial_\mu H(\widetilde{Y}_s, \mu_s)\left(Y_s\right) d \widetilde{K}_s\right]\right] \\
			\leq & \e\left[\int_0^T e^{\alpha s}\left(H(Y_s, \mu_s)-H(\widehat{y}, \delta_{\widehat{y}})\right) d K_s\right] + C \e\left[ K_T \right] \sup _{0 \leq t \leq T}\e  \left[e^{\alpha t} |\Delta Y_t| \right]
			\\
			\leq & 
			C \e\left[ |\xi|^2 + \int_0^T |f^0(s)|^2ds + |\widehat{y}|^2 + (K_T)^2 \right]. \\
		\end{aligned}
	\end{equation}
	Combining (\ref{e16}) and 
	the last inequality, from BDG inequality and Young's inequality, we obtain the desired result.
\end{proof}

\section{Stability and uniqueness of the solution}\label{sec-unique-stability}

In this section, we investigate the stability and uniqueness of the solution of (\ref{MF-RBSDE})-(\ref{refl}). First, we prove the stability of the solution in Section \ref{subsec-stability}. This stability result plays a crucial role in Section \ref{sec:obstacle problem} , as it facilitates the demonstration of the continuity of the decoupling function for the corresponding obstacle problem in Wasserstein space. 
Subsequently, we leverage the stability result to establish the uniqueness of the solution in Section \ref{subsec-uniqueness}.

\subsection{Stability result}\label{subsec-stability}

\begin{pr}\label{pr:stability}
	Assume that $(Y^i,Z^i,K^i)_{i=1,2}$ are two solutions of the mean-field reflected BSDE (\ref{MF-RBSDE})-(\ref{refl}) with data $(f^i, \xi^i)_{i=1,2}$. 
	Let Assumptions \ref{H}, \ref{H:SDE}, and \eqref{H-uniqueness} hold true for $(f^i, \xi^i)_{i=1,2}$ and $H$. Then, there exists a constant $C>0$ depending on $T$ and $L$, such that
	\begin{equation}\label{eq:est-stability}
		\begin{aligned}
			\sup_{0 \leq t \leq T} \e \left| \Delta Y_t \right|^2 
			+ \int_0^T \left| \Delta Z_s \right|^2 \,ds 
			\leq 
			C I_\delta^2, 
		\end{aligned}
	\end{equation}
	and
	\begin{equation}\label{eq:est-R}
		\begin{aligned}
			\sup_{0 \leq t \leq T} \e \left| \Delta R_t \right|
			\leq 
			C I_\delta, 
		\end{aligned}
	\end{equation}
	where $\Delta Y:= Y^1-Y^2, \Delta Z:= Z^1-Z^2, \Delta R:=R^1-R^2, \Delta \xi:= \xi^1-\xi^2, \delta f:= f^1-f^2$,
	\begin{equation}\label{eq:R-def-1}
	\begin{aligned}
		& R^i_t := \int_0^t \partial_y H(Y^i_{s}, \mu^i_{s}) \,dK^i_s 
		+ \widetilde{\e} \left[ \int_0^t  \partial_\mu H(\widetilde{Y}^i_{s}, \mu^i_{s})(Y^i_{s}) \, d \widetilde{K}^i_s \right], \quad i=1,2, 
	\end{aligned}
	\end{equation}
	and 
	\begin{equation}\label{eq:I-delta}
		\begin{aligned}
			I_\delta := \e^{\frac{1}{2} } \left[ \left| \Delta \xi \right|^2 + \int_0^T \left| \delta f(s,Y^1_s,Z^1_s,\nu^1_s) \right|^2 \,ds \right] . 
		\end{aligned}
	\end{equation}
\end{pr}

\begin{proof}
	We apply \ito to $ \left| \Delta Y_t \right|^2$.  Following the same technique used in the proof of Lemma \ref{lem-a-priori-supE}, from the concavity of $H$, we obtain 
	\begin{equation}\label{ito-delta-y}
		\begin{aligned}
			\e \left[ \left| \Delta Y_t \right|^2 + \int_t^T \left| \Delta Z_s \right|^2 \,ds \right]
			\leq 
			2 \e \left[ \int_t^T  \left| \Delta Y_s \right|  \left| \Delta f(s) \right| \,ds \right] , 
		\end{aligned}
	\end{equation}
	where $\Delta f(s):= f^1(s,Y^1_s,Z^1_s,\nu^1_s) - f^2(s,Y^2_s,Z^2_s,\nu^2_s)$. 
	
	From the Lipschitz continuity of $f^2$, we deduce that 
	\begin{equation}\label{delta-y-f}
		\begin{aligned}
			& 2 \left| \Delta Y_s \right|  \left| \Delta f(s) \right| \\
			\leq &
			2 \left| \Delta Y_s \right| 
			\left( 
			\left| \delta f(s,Y^1_s,Z^1_s,\nu^1_s) \right|
			+ \left| f^2(s,Y^1_s,Z^1_s,\nu^1_s) - f^2(s,Y^2_s,Z^2_s,\nu^2_s) \right|
			\right) \\
			\leq &
			\left| \delta f(s,Y^1_s,Z^1_s,\nu^1_s) \right|^2
			+ C_L \left| \Delta Y_s \right|^2 
			+ \e \left| \Delta Y_s \right|^2 
			+ \frac{1}{4} \left( \left| \Delta Z_s \right|^2 + \e \left| \Delta Z_s \right|^2 \right) .  \\
		\end{aligned}
	\end{equation}
	Substituting \eqref{delta-y-f} into \eqref{ito-delta-y} and applying Gronwall's lemma, we obtain \eqref{eq:est-stability}. 
	Then, \eqref{eq:est-R} is obtained from \eqref{eq:est-stability}, the equation  
	\begin{equation}\label{eq:Delta R}
		\begin{aligned}
			\Delta R_t = \Delta Y_0 - \Delta Y_t - \int_0^t \Delta f(s) \,ds + \int_0^t \Delta Z_s \,dB_s , \quad 0\leq t\leq T, 
		\end{aligned}
	\end{equation}
	and the Lipschitz continuity of $f^2$. 
\end{proof}

\subsection{Uniqueness of the solution}\label{subsec-uniqueness}

In this section, we consider the uniqueness of the solution of (\ref{MF-RBSDE})-(\ref{refl}). 

\begin{pr}\label{pr-unique}
	Let Assumption \ref{H} be satisfied. 
	Assume that $(Y^i,Z^i,K^i)_{i=1,2}$ are two solutions of the mean-field reflected BSDE (\ref{MF-RBSDE})-(\ref{refl}). Then $(Y^1,Z^1)=(Y^2,Z^2)$. 

	Moreover, if there exists a constant $\delta_0>0$, such that 
	\begin{equation}\label{H-uniqueness}
		\begin{aligned}
			\left|\e_{ X\sim\mu }\left[\partial_\mu H(y,\mu)(X)\right]\right| \leq (1-\delta_0) \left|\partial_y H(y,\mu)\right|
			,  \quad \forall(y,\mu)\in \R^n\times \Pc_2(\R^n). 
		\end{aligned}
	\end{equation}
	Then $(Y^1,Z^1,K^1)=(Y^2,Z^2,K^2)$.
\end{pr}

\begin{proof}
We set $\Delta Y:=Y^1-Y^2, \Delta Z:=Z^1-Z^2$, and $\Delta K:=K^1-K^2$. From Proposition \ref{pr:stability}, we obtain
	\begin{equation}\label{ee}
		\begin{aligned}
			\e \left[|\Delta Y_t|^2\right]
			+\frac{1}{2}\e \left[\int_t^T |\Delta Z_s|^2 ds\right]
			=
			0.
			\\
		\end{aligned}
	\end{equation}
	Hence $(Y^1,Z^1)=(Y^2,Z^2)=(Y,Z)$.

Assume moreover that \eqref{H-uniqueness} holds true. We shall prove the uniqueness of $K$. We derive from the uniqueness of $Y$ and $Z$, that for all $0\leq s\leq t\leq T$,
	\begin{equation}\label{e18}
		\begin{aligned}
			\int_s^t \partial_y H(Y_{u},\mu^{}_{u}) d\Delta K_u +\widetilde{\e}\left[\int_s^t \partial_\mu H(\widetilde{Y}_{u},\mu^{}_{u})(Y^{}_{u}) d \widetilde{\Delta K}_u\right]=0.
		\end{aligned}
	\end{equation}

	From (\ref{H-beta}) and the last equation, we obtain
	\begin{equation}\label{e2}
		\begin{aligned}
			d\Delta K_t=-\frac{\partial_y H(Y_{t},\mu^{}_{t})^*}{\left| \partial_y H(Y_{t},\mu^{}_{t})  \right|^2} \widetilde{\e}\left[ \partial_\mu H(\widetilde{Y}_{t},\mu^{}_{t})(Y^{}_{t}) d \widetilde{\Delta K}_t \right].
		\end{aligned}
	\end{equation}
Multiplying both sides by $\partial_\mu H(Y_{t},\mu^{}_{t})(\widetilde{Y}^{}_{t})$ and taking the expectation $\e [\cdot]$, we obtain
	\begin{equation}\label{e19}
	\begin{aligned}
		F(\widetilde{Y}_t) +\e\left[ G(Y_t, \widetilde{Y}_t) F(Y_t)\right] =0,
	\end{aligned}
	\end{equation}
	where 
	\begin{equation}
	\begin{aligned}
		F(Y)=
		\widetilde{\e}\left[ \partial_\mu H(\widetilde{Y},\mu)(Y) d \widetilde{\Delta K} \right]
		, \quad 
		G(Y,\widetilde{Y})=\frac{\partial_\mu H(Y,\mu)(\widetilde{Y}) \partial_y H(Y,\mu)^*}{\left| \partial_y H(Y,\mu)  \right|^2}. 
	\end{aligned}
	\end{equation}
	From (\ref{H-uniqueness}), we deduce that $\widetilde{\e}\|G(Y,\widetilde{Y})\|_2\leq 1-\delta_0$. Then, from Fubini's theorem, we obtain 
	\begin{equation}
	\begin{aligned}
		\e| F(Y_t) |
	=
		\widetilde{\e}|F(\widetilde{Y}_t)| 
	\leq 
		\e \left[\widetilde{\e}\| G(Y_t, \widetilde{Y}_t)\|_2 \left| F(Y_t) \right|\right]
	\leq (1-\delta_0) \e| F(Y_t) | . 
	\end{aligned}
	\end{equation}
	Hence, we obtain from \eqref{e2} 
	\begin{equation}
		\begin{aligned}
			\partial_y H(Y_{t},\mu^{}_{t})d\Delta K_t
			=
			\widetilde{\e}\left[ \partial_\mu H(\widetilde{Y}_{t},\mu^{}_{t})(Y^{}_{t})  d \widetilde{\Delta K}_t\right]=0.  
		\end{aligned}
	\end{equation}
Then, the uniqueness of $K$ follows from \eqref{H-beta}. 
\end{proof}

\begin{rem}
	\begin{enumerate}
		\item As discussed in Section \ref{sec-limit}, the solution of (\ref{MF-RBSDE})-(\ref{refl}) is obtained as the limit of solutions for a series of particle systems with oblique reflection. The assumption \eqref{H-uniqueness} naturally arises in the literature on obliquely reflected BSDEs. It ensures that the perturbing operator for the reflection direction is (uniformly) positive definite, satisfying \cite[Assumption (SB).(iii)]{CR20}. We refer the reader to \cite{R02,CR20} for the assumptions and well-posedness results related to obliquely reflected BSDEs.
		\item As noted in \cite[Section 1, page 484]{BEH18}, in the case where $H$ depends only on $\mu$, the reflecting process $K$ is generally not unique. Nevertheless, under the assumption that 
		\begin{equation}
		\begin{aligned}
			\inf_{\mu\in\Pc_2(\R^n)}\e_{X\sim\mu}\left| \partial_\mu H(\mu)(X) \right|^2>0, 
		\end{aligned}
		\end{equation}
		the mean process $\e[K]$ is unique \cite[Theorem 25]{BCDH20}. 
	\end{enumerate}
\end{rem}

At the end of the section, we provide a counterexample where the uniqueness does not hold if inequalities \eqref{Hx-Hy} and \eqref{H-uniqueness} are not satisfied.

\paragraph*{Counterexample.}  

Consider $n=1, \xi=1, f=0$ and $H(x,\mu)=x-\int_{\Omega} y\,\mu(dy)$. Then, it is straightforward to verify that, for all $\lambda\geq 0$, $(1,0,\lambda t)_{0\leq t\leq T}$ is a solution of (\ref{MF-RBSDE})-(\ref{refl}). Therefore, the process $K$ is not unique.

\section{Existence of the solution}\label{sec-exist}

In this section, we consider the existence of the solution of (\ref{MF-RBSDE})-(\ref{refl}). We construct a solution through a penalized BSDE approach. For $m\geq1$, let $(Y^m,Z^m)$ be the solution of the following BSDE: for all $t\in[0,T]$, 
\begin{equation}\label{penal}
	\begin{aligned}
		Y^m_t= & \xi+\int_t^T f(s,Y^m_s,Z^m_s,\nu^m_s)ds - \int_t^T Z^m_s dB_s \\
		       &
		+ \int_t^T \partial_y H(Y^m_s, \mu^m_s)dK^m_s +\widetilde{\e}\left[\int_t^T \partial_\mu H(\widetilde{Y}^m_s,\mu^m_s)(Y^m_s)d\widetilde{K}^m_s\right],
	\end{aligned}
\end{equation}
where
\begin{equation}\label{def-Km}
	\begin{aligned}
		K^m_t:=\int_0^t m H^-(Y^m_s, \mu^m_s)ds,
		\quad \mu^m_t:=[Y^m_t],\quad \nu^m_t:=[(Y^m_t,Z^m_t)], \\
	\end{aligned}
\end{equation}
and $H^-$ denotes the negative part of $H$. We then have the following (uniform) a priori estimates of $(Y^m,Z^m)$. 
\begin{lem}\label{lem-a priori-supE-Ym}
	Suppose that Assumption \ref{H} holds true. Then, for any $m\geq 1$, the penalized equation (\ref{penal}) admits a unique solution $(Y^m,Z^m)\in S^{2,n}\times H^{2,n}$. Moreover, there exists a constant $C>0$ independent of $m$, such that
	\begin{equation}\label{est-supE-Ym}
		\begin{aligned}
			\sup_{0\leq t \leq T} \e\left[	\left| Y^m_t \right|^2\right] + \e\left[ \int_0^T \left|Z^m_s\right|^2ds \right]\leq C \e \left[ |\xi|^2 + \int_0^T |f^0(s)|^2ds + |\widehat{y}|^2 \right],
		\end{aligned}
	\end{equation}
	\begin{equation}\label{est-supE-Ym4}
		\begin{aligned}
			\sup_{0\leq t \leq T}\e\left[ 	\left| Y^m_t \right|^4\right] \leq C \e \left[ |\xi|^4 + \int_0^T |f^0(s)|^4ds + |\widehat{y}|^4 \right].
		\end{aligned}
	\end{equation}
\end{lem}

\begin{proof}
	The existence and uniqueness of a solution $(Y^m,Z^m)$ of (\ref{penal}) follow from a (slightly) generalized version of \cite[Theorem 3.1]{BLP09}. 

	It can easily be deduced from the definition of $K^m$ that for all $t\in[0,T]$,
	\begin{equation}
		\begin{aligned}
			\int_t^T H(Y^m_s, \mu^m_s)dK^m_s\leq0,
		\end{aligned}
	\end{equation}
	and the estimates (\ref{est-supE-Ym}) and \eqref{est-supE-Ym4} are obtained in an identical way to the proof of Lemma \ref{lem-a-priori-supE}.
\end{proof}

\begin{lem}\label{lem-eHm}
	Suppose that Assumption \ref{H} holds true. Then, there exists a constant $C>0$, such that for all $m\geq 1$,
	\begin{equation}\label{est-H^-}
		\begin{aligned}
			\e \left[\sup_{0\leq t \leq T} H^-(Y^m_t,\mu^m_t)^2\right]\leq
			\frac{C}{m} \,\,\,\mathrm{and}\,\,\,
			\e\left[\int_0^T  H^-(Y^m_s,\mu^m_s)^2ds\right]
			\leq
			\frac{C}{m^2}.
		\end{aligned}
	\end{equation}
	In particular, for all $m\geq1$,
	\begin{equation}\label{EKm}
		\begin{aligned}
			\e \left(K^m_T\right)^2\leq C.
		\end{aligned}
	\end{equation}
\end{lem}

\begin{proof}
	The main idea of the proof is applying \ito to $H^-(Y^m_t,\mu^m_t)^2$ and using the concavity of $H$ and \eqref{est-supE-Ym} to obtain the desired estimate. However, due to the lack of needed integrability of the process $Z$, we cannot directly apply \ito to $H^-(Y^m_t,\mu^m_t)^2$. Instead, we replace $\mu_t$ with the empirical distribution, and apply the classical Itô's formula. Define $(Y^{m,i},Z^{m,i})_{1\leq i\leq N}$ as $N$ i.i.d. copies of $(Y^m,Z^m)$. Let $\mu^{m,N}_t:=N^{-1}\sum_{i=1}^N \delta_{Y^{m,i}_t}$. Applying \ito to $H^-(Y^m_t,\mu^{m,N}_t)^2$ and using the concavity of $H$, we obtain 
	\begin{equation}\label{s4l2-e1}
		\begin{aligned}
			  & H^-(Y^m_t,\mu^{m,N}_t)^2-H^-(\xi,\mu^{m,N}_T)^2
			+\int_t^T
			\mathbf{1}_{\{ H(Y^m_s,\mu^{m,N}_s)\leq0 \}}
			\left|\partial_y H(Y^m_s,\mu^{m,N}_s)^* Z^m_s\right|^2ds                                                                               \\
			  & +\frac{1}{N^2}\sum_{i=1}^N \int_t^T
			\mathbf{1}_{\{ H(Y^m_s,\mu^{m,N}_s)\leq0 \}}\left|\partial_\mu H(Y^m_s,\mu^{m,N}_s)(Y^{m,i}_s)^* Z^{m,i}_s\right|^2ds
			\\
		\leq & 
			\sum_{i=1}^{2} \Ec^f_i + \sum_{i=1}^{2}  \Ec^\Mc_i + \sum_{i=1}^{4} \Ec^R_i, 
		\end{aligned}
	\end{equation}
	where 
	\begin{equation}\label{eq:}
	\begin{aligned}
		\Ec^f_1 &= -\int_t^T 2 H^-(Y^m_s,\mu^{m,N}_s)  \partial_y H(Y^m_s,\mu^{m,N}_s) \cdot f(s,Y^m_s,Z^m_s,\nu^m_s)ds, \\
		\Ec^f_2 &= -\frac{1}{N}\sum_{i=1}^N \int_t^T 2 H^-(Y^m_s,\mu^{m,N}_s) \partial_\mu H(Y^m_s,\mu^{m,N}_s)(Y^{m,i}_s) \cdot f^i(s,Y^{m,i}_s,Z^{m,i}_s,\nu^m_s) ds, \\
		\Ec^\Mc_1 &= \int_t^T 2 H^-(Y^m_s,\mu^{m,N}_s) \partial_y H(Y^m_s,\mu^{m,N}_s)\cdot Z^m_sdB_s, \\
		\Ec^\Mc_2 &= \frac{1}{N}\sum_{i=1}^N\int_t^T 2 H^-(Y^m_s,\mu^{m,N}_s) \partial_\mu H(Y^m_s,\mu^{m,N}_s)(Y^{m,i}_s)\cdot Z^{m,i}_sdB^i_s, 
		\\
		\Ec^R_1 &= -\int_t^T 2 H^-(Y^m_s,\mu^{m,N}_s)  \partial_y H(Y^m_s,\mu^{m,N}_s) \cdot \partial_y H(Y^m_s,\mu^{m}_s) dK^m_s, \\
		\Ec^R_2 &= -\int_t^T 2 H^-(Y^m_s,\mu^{m,N}_s)  \partial_y H(Y^m_s,\mu^{m,N}_s) \cdot \widetilde{\e}\left[\partial_\mu H(\widetilde{Y}^m_s,\mu^{m}_s)(Y^m_s) d\widetilde{K}^m_s\right], \\
		\Ec^R_3 &= -\frac{1}{N}\sum_{i=1}^N
		\int_t^T 2 H^-(Y^m_s,\mu^{m,N}_s)
		\partial_\mu H(Y^{m}_s,\mu^{m,N}_s)(Y^{m,i}_s)\cdot
		\partial_y H(Y^{m,i}_s,\mu^{m}_s)dK^i_s, \\
		\Ec^R_4 &= -\frac{1}{N}\sum_{i=1}^N
		\int_t^T 2 H^-(Y^m_s,\mu^{m,N}_s)
		\partial_\mu H(Y^{m}_s,\mu^{m,N}_s)(Y^{m,i}_s)\cdot
		\widetilde{\e}\left[\partial_\mu H(\widetilde{Y}^m_s,\mu^{m}_s)(Y^{m,i}_s) d\widetilde{K}^m_s\right]. \\
	\end{aligned}
	\end{equation}

	We now study each term seperately. 
	
	We start by dealing with the reflection terms $\Ec^R:=\sum_{1\leq i\leq 4} \Ec^R_i$. From \eqref{H-beta} and \eqref{def-Km}, we obtain
	\begin{equation}\label{s4l2-e1-k}
		\begin{aligned}
			\Ec^R \leq & \Ec^R_1 \\
			= & -\int_t^T 2 H^-(Y^m_s,\mu^{m}_s)  \left|\partial_y H(Y^m_s,\mu^{m}_s)\right|^2  dK^m_s \\
			& + \int_t^T 2 \left( H^-(Y^m_s,\mu^{m}_s) \partial_y H(Y^m_s,\mu^{m}_s) - H^-(Y^m_s,\mu^{m,N}_s) \partial_y H(Y^m_s,\mu^{m,N}_s) \right) \cdot \partial_y H(Y^m_s,\mu^{m}_s) dK^m_s
			\\
			\leq &  -2\beta^2 m \int_t^T H^-(Y^m_s,\mu^{m}_s)^2 ds \\
			& + 2 L M m \int_t^T H^-(Y^m_s,\mu^{m}_s) \left( H^-(Y^m_s,\mu^{m}_s) + \left| \partial_y H(Y^m_s,\mu^{m,N}_s) \right| \right) W_2(\mu^{m,N}_s,\mu^m_s) ds. \\
		\end{aligned}
	\end{equation}

	Next, the terms $\Ec^f:=\sum_{1\leq i\leq 2}\Ec^f_i$ are handled using the same technique as in the proof of \cite[Lemma 5.4]{GPP96}. From Young's inequality, we have 
	\begin{equation}\label{s4l2-e1-f}
		\begin{aligned}
			\Ec^f
			\leq & 2M\int_t^T H^-(Y^m_s,\mu^{m,N}_s)\left( |f(s,Y^m_s,Z^m_s,\nu^m_s)|+\frac{1}{N}\sum_{i=1}^N |f(s,Y^{m,i}_s,Z^{m,i}_s,\nu^m_s)| \right)ds
			\\
			\leq &
			\frac{C}{m} \int_t^T \left( |f(s,Y^m_s,Z^m_s,\nu^m_s)|^2 + \frac{1}{N}\sum_{i=1}^N |f^i(s,Y^{m,i}_s,Z^{m,i}_s,\nu^m_s)|^2 \right)ds
			\\
			     & + \beta^2 m \int_t^T H^-(Y^m_s,\mu^{m,N}_s)^2ds . 
		\end{aligned}
	\end{equation}

	Bring together the estimates on $\Ec^R$ and $\Ec^f$, we obtain 
	\begin{equation}\label{s4l2-e2}
		\begin{aligned}
			     & H^-(Y^m_t,\mu^{m,N}_t)^2-H^-(\xi,\mu^{m,N}_T)^2
			+\int_t^T
			\mathbf{1}_{\{ H(Y^m_s,\mu^{m,N}_s)\leq0 \}}
			\left|\partial_y H(Y^m_s,\mu^{m,N}_s)^* Z^m_s\right|^2ds                                                         \\
			     & +\frac{1}{N^2}\sum_{i=1}^N \int_t^T
			\mathbf{1}_{\{ H(Y^m_s,\mu^{m,N}_s)\leq0 \}}\left|\partial_\mu H(Y^m_s,\mu^{m,N}_s)(Y^{m,i}_s)^* Z^{m,i}_s\right|^2ds
			\\
			\leq &
			\frac{C}{m} \int_t^T \left( |f(s,Y^m_s,Z^m_s,\nu^m_s)|^2 + \frac{1}{N}\sum_{i=1}^N |f^i(s,Y^{m,i}_s,Z^{m,i}_s,\nu^m_s)|^2 \right)ds \\
			& + 2 L M m \int_t^T W_2(\mu^{m,N}_s,\mu^m_s) H^-(Y^m_s,\mu^{m}_s) \left( H^-(Y^m_s,\mu^{m}_s) + M \right)ds \\
			& -\beta^2 m \int_t^T H^-(Y^m_s,\mu^{m}_s)^2 ds 
			+\Ec^\Mc, 
			\\
		\end{aligned}
	\end{equation}
	where $\Ec^\Mc=\Ec^\Mc_1+\Ec^\Mc_2$. From the Lipschitz continuity of $H$ and Lemma \ref{lem-a priori-supE-Ym}, we deduce that the local martingale $\Ec^\Mc$ is a true martingale. 

	Then, we take the expectation in \eqref{s4l2-e2} and let $N\to\infty$. From \eqref{H-M}, \eqref{Hy-L} and \eqref{est-supE-Ym}, we have 
	\begin{equation}\label{s4l2-e10}
	\begin{aligned}
		& \e\left|H^-(Y^m_t,\mu^{m,N}_t)^2-H^-(Y^m_t,\mu^{m}_t)^2\right| \\
	\leq &
		\e^{\frac{1}{2}}\left|H^-(Y^m_t,\mu^{m,N}_t)+H^-(Y^m_t,\mu^{m}_t)\right|^2
		\e^{\frac{1}{2}}\left|H^-(Y^m_t,\mu^{m,N}_t)-H^-(Y^m_t,\mu^{m}_t)\right|^2 \\
	\leq &
		C \e^{\frac{1}{2}} \left[ \left( |Y^m_t|+ \frac{1}{N}\sum_{i=1}^N |Y^{m,i}_t| \right)^2 \right] 
		\e^{\frac{1}{2}}\left[ W^2_2(\mu^{m,N}_t,\mu^m_t) \right] \\
	\leq &
		C \e^{\frac{1}{2}}\left[ W^2_2(\mu^{m,N}_t,\mu^m_t) \right], \\
	\end{aligned}
	\end{equation}
	for some constant $C\geq0$ independent of $N$. 

	From Lipschitz continuity of $H$ and \eqref{est-supE-Ym4}, we have 
	\begin{equation}\label{s4l2-e11}
	\begin{aligned}
		& \e\left[\int_0^T W_2(\mu^{m,N}_s,\mu^m_s) H^-(Y^m_s,\mu^{m}_s)^2ds\right] \\
	\leq & 
		C \e^{\frac{1}{2}}\left[\int_0^T W^2_2(\mu^{m,N}_s,\mu^m_s) ds \right] 
		\e^{\frac{1}{2}}\left[ H^-(Y^m_s,\mu^{m}_s)^4 \right] \\
	\leq & 
		C \e^{\frac{1}{2}}\left[\int_0^T W^2_2(\mu^{m,N}_s,\mu^m_s) ds \right]. 
	\end{aligned}
	\end{equation}

	Hence, passing to the limit in \eqref{s4l2-e2}, from \eqref{s4l2-e10}, \eqref{s4l2-e11} and the limit
	\begin{equation}\label{eq:limit}
	\begin{aligned}
		\e\left[\int_0^T W^2_2(\mu^{m,N}_s,\mu^m_s) ds\right]\rightarrow0,
	\end{aligned}
	\end{equation}
	we obtain
	\begin{equation}\label{s4l2-e3}
		\begin{aligned}
			     & H^-(Y^m_t,\mu^{m}_t)^2
			+\int_t^T
			\mathbf{1}_{\{ H(Y^m_s,\mu^{m}_s)\leq0 \}}
			\left|\partial_y H(Y^m_s,\mu^{m}_s)^* Z^m_s\right|^2ds                                             \\
			\leq &
			\frac{C}{m} \int_t^T \left( |f(s,Y^m_s,Z^m_s,\nu^m_s)|^2 + \e\left[ |f(s,Y^{m}_s,Z^{m}_s,\nu^m_s)|^2\right] \right)ds \\ 
			&+\int_t^T 2 H^-(Y^m_s,\mu^{m}_s) \partial_y H(Y^m_s,\mu^{m}_s)\cdot Z^m_sdB_s 		
			-\beta^2 m \int_t^T H^-(Y^m_s,\mu^{m}_s)^2ds .  \\
		\end{aligned}
	\end{equation}

	Taking the expectation, we have 
	\begin{equation}\label{s4l2-e4}
		\begin{aligned}
			\sup_{ 0 \leq t \leq T } \e \left[  H^-(Y^m_t,\mu^{m}_t)^2 \right]
			+\e\left[\int_0^T
			\mathbf{1}_{\{ H(Y^m_s,\mu^{m}_s)\leq0 \}}
			\left|\partial_y H(Y^m_s,\mu^{m}_s)^* Z^m_s\right|^2ds     \right] 
			\leq 
			\frac{C}{m}, 
			\\
		\end{aligned}
	\end{equation}
	and 
	\begin{equation}\label{s4l2-e5}
	\begin{aligned}
		\e \left[ \int_0^T H^-(Y^m_s,\mu^{m}_s)^2ds \right] 
		\leq \frac{C}{m^2}. 
		\\ 
	\end{aligned}
	\end{equation}

	Coming back to (\ref{s4l2-e3}), we take the supremum in time and the expectation. From BDG inequality and Young's inequality, we obtain 
	\begin{equation}\label{e13-B}
		\begin{aligned}
			& \e\left[\sup_{0\leq t\leq T}\left| 
			\int_t^T H^-(Y^m_s,\mu^{m}_s) \partial_y H(Y^m_s,\mu^{m}_s)\cdot Z^m_sdB_s
			\right| \right] \\
			\leq &
			\e\left[ \int_0^T \mathbf{1}_{\{H(Y^m_s,\mu^m_s)\leq0\}} |H^-(Y^m_s,\mu^{m}_s)|^2 |\partial_y H(Y^m_s,\mu^{m}_s)^* Z^m_s|^2 ds \right]^{\frac{1}{2}} \\
			\leq & 
			\epsilon \e\left[\sup_{0\leq t \leq T} |H^-(Y^m_t,\mu^m_t)|^2 \right]
			+ C_\epsilon \e\left[\int_0^T \mathbf{1}_{\{H(Y^m_s,\mu^{m}_s)\leq0\}} |H_1(Y^m_s,\mu^{m}_s)^* Z^m_s|^2 ds\right] 
			. 
		\end{aligned}
	\end{equation}
	
	Then, from (\ref{s4l2-e3}),  (\ref{s4l2-e4}) and the last inequality, we obtain  
	\begin{equation}\label{e13-Esup} 
		\begin{aligned}
			& \e\left[\sup_{0\leq t\leq T}H^-(Y^m_t,\mu^{m}_t)^2\right]+\e\left[\int_0^T \mathbf{1}_{\{H(Y^m_s,\mu^{m}_s)\leq0\}} |\partial_y H(Y^{m}_{s},\mu^{m}_{s}) ^* Z^m_s|^2ds\right] 
			\leq  
			\frac{C}{m}, 
		\end{aligned}
	\end{equation}
	which, together with \eqref{s4l2-e5}, proves \eqref{est-H^-}. 
\end{proof}

\begin{lem}\label{lem-a priori-Esup Ym}
	Under Assumption \ref{H}, there exists a constant $C>0$ independent of $m$, such that
	\begin{equation}
		\begin{aligned}
			\e\left[\sup_{0\leq t \leq T}	\left| Y^m_t \right|^2 + \int_0^T \left|Z^m_s\right|^2ds \right]\leq C \e \left[ |\xi|^2 + \int_0^T |f^0(s)|^2ds + |\widehat{y}|^2 + 1 \right].
		\end{aligned}
	\end{equation}
\end{lem}

\begin{proof}
	The result can be obtained in an identical way to the proof of Lemma \ref{lem-a-priori-Esup}. 
	Therefore, we omit the proof. 
\end{proof}

\begin{lem}\label{lem-Ym-Zm-Cauchy}
	Under Assumption \ref{H}, the sequence $(Y^m,Z^m)_{m\geq 1}$ is a Cauchy sequence in $S^{2,n} \times H^{2,n}$.
\end{lem}

\begin{proof}
	The proof follows mainly from arguments in the proof of \cite[Lemma 5.6]{GPP96}. However, some extra work is required to handle the reflection terms properly. 
	Let $m,l$ be two positive integers, and set $\Delta Y := Y^m -Y^l$ and $\Delta Z := Z^m -Z^l$. We apply \ito to $e^{\alpha t}|\Delta Y_t|^2$ and obtain
	\begin{equation}\label{ito-Ym-Yl}
		\begin{aligned}
			  & e^{\alpha t}|\Delta Y_t|^2 + \int_t^T  e^{\alpha s}|\Delta Z_s|^2ds \\
			= &
			- \int_t^T \alpha e^{\alpha s} |\Delta Y_s|^2ds+ 2\int_t^T e^{\alpha s}\Delta Y_s \cdot \Delta f(s)ds - 2 \int_t^T
			e^{\alpha s}\Delta Y_s \cdot \Delta Z_s
			dB_s                                                                    \\
			  & + 2  \int_t^T
			e^{\alpha s}\Delta Y_s \cdot
			\left( \partial_y H(Y^m_s,\mu^m_s)dK^m_s - \partial_y H(Y^l_s,\mu^l_s)dK^l_s  \right)
			\\
			  & + 2 \int_t^T
			e^{\alpha s}\Delta Y_s \cdot
			\left( \widetilde{\e}\left[\partial_\mu H(\widetilde{Y}^m_s,\mu^m_s)(Y^m_s)d\widetilde{K}^m_s\right] - \widetilde{\e}\left[\partial_\mu H(\widetilde{Y}^l_s,\mu^l_s)(Y^l_s)d\widetilde{K}^l_s\right]  \right) ,
			\\
		\end{aligned}
	\end{equation}
	where $\Delta f(s)=f(s,Y^m_s,Z^m_s,\nu^m_s)-f(s,Y^l_s,Z^l_s,\nu^l_s)$.

	Taking the expectation on both sides of (\ref{ito-Ym-Yl}), we have
	\begin{equation}\label{eito-Ym-Yl}
		\begin{aligned}
			     & \e \left[	e^{\alpha t} |\Delta Y_t|^2 + \frac{1}{2}\int_t^T  e^{\alpha s} |\Delta Z_s|^2ds \right] \\
			\leq &
			2 \e \left[ \int_t^T
				e^{\alpha s}\Delta Y_s \cdot
				\left( \partial_y H(Y^m_s,\mu^m_s)dK^m_s - \partial_y H(Y^l_s,\mu^l_s)dK^l_s  \right)
				\right]
			\\
			     & + 2\e\left[ \int_t^T
				e^{\alpha s}\Delta Y_s \cdot
				\left( \widetilde{\e}\left[\partial_\mu H(\widetilde{Y}^m_s,\mu^m_s)(Y^m_s)d\widetilde{K}^m_s\right] - \widetilde{\e}\left[\partial_\mu H(\widetilde{Y}^l_s,\mu^l_s)(Y^l_s)d\widetilde{K}^l_s\right]  \right)\right].
			\\
		\end{aligned}
	\end{equation}

	From the concavity of $H$, Lemma \ref{lem-eHm} and Fubini's theorem, we obtain 
	\begin{equation}\label{est-Km-Kl-1}
		\begin{aligned}
			     & \e \left[\int_t^T e^{\alpha s}\Delta Y_s \cdot
				\left(
				\partial_y H(Y^m_s,\mu^m_s)dK^m_s
				+
				\widetilde{\e}\left[\partial_\mu H(\widetilde{Y}^m_s,\mu^m_s)(Y^m_s)d\widetilde{K}^m_s\right]
				\right)
				\right]
			\\
			=    &
			\e \left[ \int_t^T
				e^{\alpha s}\Delta Y_s\cdot
				\partial_y H(Y^m_s,\mu^m_s)dK^m_s
				+ \int_t^T e^{\alpha s}
				\widetilde{\e}\left[ \Delta\widetilde{Y}_s\cdot \partial_\mu H(Y^m_s,\mu^m_s)(\widetilde{Y}^m_s)\right]dK^m_s
			\right]                                                                                                                       \\
			&  + \e \left[ 
			\int_t^T e^{\alpha s}\Delta Y_s\cdot\widetilde{\e}\left[\partial_\mu H(\widetilde{Y}^m_s,\mu^m_s)(Y^m_s)d\widetilde{K}^m_s\right]	
				 - \int_t^T
			e^{\alpha s}\widetilde{\e}\left[ \Delta\widetilde{Y}_s\cdot \partial_\mu H(Y^m_s,\mu^m_s)(\widetilde{Y}^m_s)\right]dK^m_s
			\right]
			\\
			\leq &
			\e \left[\int_t^T e^{\alpha s}
			\left( H(Y^m_s,\mu^m_s) - H(Y^l_s,\mu^l_s) \right) dK^m_s\right]                                                        \\
			\leq &
			\e \left[\int_t^T e^{\alpha s}
			H^-(Y^l_s,\mu^l_s) dK^m_s\right]                                                                                          \\
			=    &
			m \e \left[\int_t^T e^{\alpha s}
			H^-(Y^l_s,\mu^l_s)H^-(Y^m_s,\mu^m_s) ds\right]                                                                           \\
			\leq &
			e^{\alpha T} m \,  \e^{\frac{1}{2} } \left[\int_t^T H^-(Y^l_s,\mu^l_s)^2ds \right] \e^{\frac{1}{2} } \left[\int_t^T H^-(Y^m_s,\mu^m_s)^2ds \right]  \\
			\leq &
			\frac{C}{l} ,                                                                      \\
		\end{aligned}
	\end{equation}
	where $C>0$ is the constant appeared in Lemma \ref{lem-eHm} multiplied by $e^{\alpha T}$, which is independent of $m$ and $l$.

	Arguing similarly, we also have
	\begin{equation}\label{est-Km-Kl-2}
		\begin{aligned}
			 & - \e \left[\int_t^T e^{\alpha s}\Delta Y_s \cdot
				\left(
				\partial_y H(Y^l_s,\mu^l_s)dK^l_s
				+
				\widetilde{\e}\left[\partial_\mu H(\widetilde{Y}^l_s,\mu^l_s)(Y^l_s)d\widetilde{K}^l_s\right]
				\right)
				\right]
			\leq
			\frac{C}{m}.                                        \\
		\end{aligned}
	\end{equation}

	Combining (\ref{eito-Ym-Yl}), (\ref{est-Km-Kl-1}), and (\ref{est-Km-Kl-2}), we obtain
	\begin{equation}\label{est-supE-Km-Kl}
		\begin{aligned}
			 & \sup_{0\leq t\leq T} \e \left[  |\Delta Y_t|^2\right] + \e \left[ \int_0^T  |\Delta Z_s|^2ds \right]
			\leq
			C \left( \frac{1}{l} + \frac{1}{m} \right) . 
		\end{aligned}
	\end{equation}

	Coming back once again to (\ref{ito-Ym-Yl}), taking the supremum in $t$ and the expectation, we obtain
	\begin{equation}\label{e14-Esup}
		\begin{aligned}
			     & \e\left[\sup_{0\leq t \leq T}e^{\alpha t}|\Delta Y_t|^2 + \frac{1}{2}\int_0^T  e^{\alpha s}|\Delta Z_s|^2ds\right] \\
			\leq &
			2 \e\left[\sup_{0\leq t \leq T}\left|\int_t^T
				e^{\alpha s}\Delta Y_s\cdot\Delta Z_s
			dB_s\right|\right]                                                                                                        \\
			     & + 2  \e\sup_{0\leq t \leq T}\left[\int_t^T e^{\alpha s}\Delta Y_s \cdot
				\left(
				\partial_y H(Y^m_s,\mu^m_s)dK^m_s
				+
				\widetilde{\e}\left[\partial_\mu H(\widetilde{Y}^m_s,\mu^m_s)(Y^m_s)d\widetilde{K}^m_s\right]
				\right)
			\right]                                                                                                                     \\
			     & + 2  \e\sup_{0\leq t \leq T}\left[-\int_t^T e^{\alpha s}\Delta Y_s \cdot
				\left(
				\partial_y H(Y^l_s,\mu^l_s)dK^l_s
				+
				\widetilde{\e}\left[\partial_\mu H(\widetilde{Y}^l_s,\mu^l_s)(Y^l_s)d\widetilde{K}^l_s\right]
				\right)
			\right].                                                                                                                    \\
		\end{aligned}
	\end{equation}

	From BDG inequality and Young's inequality, we have
	\begin{equation}\label{e14-B}
		\begin{aligned}
			\e\left[\sup_{0\leq t \leq T}\left|\int_t^T
				e^{\alpha s}\Delta Y_s\cdot\Delta Z_s
				dB_s\right|\right]
			\leq
			\epsilon \e\left[\sup_{0\leq t \leq T} e^{\alpha t} |\Delta Y_t|^2 \right]
			+ C_\epsilon \e\left[ \int_0^T e^{\alpha s}|\Delta Z_s|^2ds \right].
		\end{aligned}
	\end{equation}

	From the concavity of $H$ and the boundedness of $\partial_y H$ and $\partial_\mu H$, thanks to Lemma \ref{lem-eHm}, we have
	\begin{equation}\label{e14-Km-Lm}
		\begin{aligned}
			     & \e\sup_{0\leq t \leq T}\left[\int_t^T e^{\alpha s}\Delta Y_s \cdot
				\left(
				\partial_y H(Y^m_s,\mu^m_s)dK^m_s
				+
				\widetilde{\e}\left[\partial_\mu H(\widetilde{Y}^m_s,\mu^m_s)(Y^m_s)d\widetilde{K}^m_s\right]
				\right)
			\right]                                                                                                                                                                          \\
			\leq &
			\e \sup_{0\leq t \leq T}\left[ \int_t^T
				e^{\alpha s} \Delta Y_s \cdot
				\partial_y H(Y^m_{s},\mu^{m}_{s}) dK^m_s
				+ \int_t^T e^{\alpha s}
				\widetilde{\e}\left[ \Delta Y_s\cdot \partial_\mu H(Y^m_s,\mu^m_s)(\widetilde{Y}^m_s)\right]dK^m_s
			\right]                                                                                                                                                                          \\
			     & + \e \sup_{0\leq t \leq T}\left[ -\int_t^T e^{\alpha s}
				\widetilde{\e}\left[ \Delta\widetilde{Y}_s\cdot \partial_\mu H(Y^m_s,\mu^m_s)(\widetilde{Y}^m_s)\right]dK^m_s
			\right]                                                                                                                                                                          \\
			     & +\e\sup_{0\leq t \leq T}\left[\int_t^T e^{\alpha s}\Delta Y_s\cdot\widetilde{\e}\left[\partial_\mu H(\widetilde{Y}^m_s,\mu^m_s)(Y^m_s)d\widetilde{K}^m_s\right] \right] \\
			\leq &
			\e \sup_{0\leq t \leq T}\left[ \int_t^T
				e^{\alpha s}
				\left( H(Y^m_s,\mu^m_s) - H(Y^l_s,\mu^l_s) \right) dK^m_s \right]
			+ 2M \sup_{0\leq t \leq T}\e \left[ e^{\alpha t}|\Delta Y_t|\right] \e\left[K^m_T\right]
			\\
			\leq &
			C \left(\frac{1}{l}+\frac{1}{\sqrt{l}}+\frac{1}{\sqrt{m}}\right).                                                                                                                \\
		\end{aligned}
	\end{equation}

	Arguing similarly, we also have
	\begin{equation}\label{e14-Kl-Ll}
		\begin{aligned}
			 & \e\sup_{0\leq t \leq T}\left[ - \int_t^T
				e^{\alpha s}\Delta Y_s \cdot \left(
				\partial_y H(Y^{l}_{s},\mu^{l}_{s}) dK^l_s + \widetilde{\e} \left[ \partial_\mu H(\widetilde{Y}^{l}_{s},\mu^{l}_{s})(Y^{l}_{s})d \widetilde{K}^l_S  \right]
				\right)
				\right]
			\\
			\leq &
			C \left(\frac{1}{m}+\frac{1}{\sqrt{l}}+\frac{1}{\sqrt{m}}\right). \\
		\end{aligned}
	\end{equation}

	Combining (\ref{e14-Esup}), (\ref{e14-B}), (\ref{e14-Km-Lm}) and (\ref{e14-Kl-Ll}), we deduce that
	\begin{equation}\label{e15}
		\begin{aligned}
			 & \e\left[\sup_{0\leq t \leq T}e^{\alpha t}|\Delta Y_t|^2 + \frac{1}{2}\int_0^T  e^{\alpha s}|\Delta Z_s|^2ds\right]
			\leq
			C \left(\frac{1}{\sqrt{l}}+\frac{1}{\sqrt{m}}\right).                                                                   \\
		\end{aligned}
	\end{equation}

	Therefore, we conclude that $(Y^m,Z^m)_{m\geq1}$ is a Cauchy sequence in $S^{2,n}\times H^{2,n}$.

\end{proof}

We now have all the key ingredients to show the existence result. 
\begin{pr}\label{pr-exist}
	Under Assumption \ref{H}, there exists a solution $(Y,Z,K)$ of (\ref{MF-RBSDE})-(\ref{refl}) in $S^{2,n}\times H^{2,n}\times A^{2,1}$.
\end{pr}

\begin{proof}
	It follows from Lemma \ref{lem-Ym-Zm-Cauchy} that
	\begin{equation}
		\begin{aligned}
			Y^m\overset{S^{2,n}}{\longrightarrow}Y, \quad Z^m\overset{H^{2,n}}{\longrightarrow}Z, \,\,\, \mathrm{and}\,\,\, \int_0^\cdot Z^m_sdB_s\overset{L^{2,n}}{\longrightarrow}\int_0^\cdot Z_sdB_s.
		\end{aligned}
	\end{equation}
	Set $\mu:=[(Y)]$ and $\nu:=[(Y,Z)]$. From the Lipschitz continuity of $f$, we have
	\begin{equation}
		\begin{aligned}
			f(\cdot,Y^m_\cdot,Z^m_\cdot,\nu^m_\cdot)\overset{H^{2,n}}{\longrightarrow}f(\cdot,Y_\cdot,Z_\cdot,\nu_\cdot),
		\end{aligned}
	\end{equation}
	Define $k^m_\cdot:=m H^-(Y^m_\cdot,\mu^m_\cdot)$. Thanks to (\ref{est-H^-}), we have $\e\left[ \int_0^T |k^m_s|^2ds \right]\leq C$. Hence, there is a subsequence of $k^m$ converging weakly to some $k$ in $L^2(\Omega\times[0,T])$. From Mazur's Lemma, we know that there exists a convex combination of $k^m$ converging strongly to $k$, namely
	\begin{equation}
		\begin{aligned}
			\sum_{i=m}^{N_m} \lambda^m_i k^i
			\overset{L^{2}}{\longrightarrow}
			k
			,
		\end{aligned}
	\end{equation}
	where $\lambda^m_i\geq0$ for all $m\geq1$ and $m\leq i\leq N_m$, and $\sum_{i=m}^{N_m}\lambda^m_i=1$.

	For all $m\geq1$, we have
	\begin{equation}
		\begin{aligned}
			\sum_{i=m}^{N_m} \lambda^m_i \partial_y H(Y^i,\mu^i)  k^i
			= \sum_{i=m}^{N_m} \lambda^m_i \left(\partial_y H(Y^i,\mu^i)-\partial_y H(Y,\mu)\right)  k^i + \partial_y H(Y,\mu) \sum_{i=m}^{N_m}\lambda^m_i  k^i.
		\end{aligned}
	\end{equation}

	From the Lipschitz continuity of $\partial_y H$, the strong $S^{2,n}$-convergence of $Y$ and the uniform $L^2(\Omega\times[0,T])$-boundedness of $k^i$, we obtain
	\begin{equation}
		\begin{aligned}
			\sum_{i=m}^{N_m} \lambda^m_i \left(\partial_y H(Y^i,\mu^i)-\partial_y H(Y,\mu)\right)  k^i
			\overset{L^{2}}{\longrightarrow}0,
		\end{aligned}
	\end{equation}
	and
	\begin{equation}
		\begin{aligned}
			\partial_y H(Y,\mu) \sum_{i=m}^{N_m}\lambda^m_i  k^i
			\overset{L^{2}}{\longrightarrow} \partial_y H(Y,\mu) k.
		\end{aligned}
	\end{equation}

	Arguing similarly, we also have
	\begin{equation}
		\begin{aligned}
			\sum_{i=m}^{N_m} \lambda^m_i \widetilde{\e}\left[\partial_\mu H(\widetilde{Y}^i,\mu^i)(Y^i)  \widetilde{k}^i\right]
			\overset{L^{2}}{\longrightarrow} \widetilde{\e}\left[\partial_\mu H(\widetilde{Y},\mu)(Y) \widetilde{k}\right].
		\end{aligned}
	\end{equation}

	Setting $K_t:=\int_0^t k_sds$, and passing to the limit into
	\begin{equation}
		\begin{aligned}
			\sum_{i=m}^{N_m} \lambda^m_i Y^i_s
			= & \xi+ \sum_{i=m}^{N_m} \int_t^T  \lambda^m_i f(s,Y^i_s,Z^i_s,\nu^i_s)ds
			- \sum_{i=m}^{N_m} \int_t^T  \lambda^m_i Z^i_sdB_s                         \\
			  & + \sum_{i=m}^{N_m} \int_t^T  \lambda^m_i \partial_y H(Y^i_s,\mu^i_s)dK^i_s
			+ \sum_{i=m}^{N_m} \e\left[ \int_t^T  \lambda^m_i \partial_\mu H(\widetilde{Y}^i_s,\mu^i_s)(Y^i_s)d\widetilde{K}^i_s \right],
		\end{aligned}
	\end{equation}
	we conclude that $(Y,Z,K)$ satisfies the BSDE (\ref{MF-RBSDE}).

	Lastly, we need to check that $H(Y_t,\mu_t)\geq0$, and that the Skorokhod condition is satisfied. Indeed, we have 
	\begin{equation}
		\begin{aligned}
			H(Y_t,\mu_t) = \lim_{m\to\infty} H(Y^m_t,\mu^m_t) \geq 0.
		\end{aligned}
	\end{equation}
	From the inequality $\int_0^T H(Y^i_s,\mu^i_s)dK^i_s\leq0$, we have
	\begin{equation}
		\begin{aligned}
			0\leq \int_0^T H(Y_s,\mu_s)dK_s = \lim_{m\to\infty}\sum_{i=m}^{N_m}
			\int_0^T \lambda^m_i H(Y^i_s,\mu^i_s)dK^i_s\leq0,
		\end{aligned}
	\end{equation}
	which verifies the Skorokhod condition.
\end{proof}

Combining Proposition \ref{pr-unique} and Proposition \ref{pr-exist}, we are ready to state the main theorem of well-posedness of (\ref{MF-RBSDE})-(\ref{refl}). 
\begin{thm}\label{thm-exist-unique}  
	Assume that Assumption \ref{H} holds true. Then, there exists a solution $(Y,Z,K)$ to the mean-field reflected BSDE (\ref{MF-RBSDE})-(\ref{refl}), and the tuple $(Y,Z)$ is unique.  If (\ref{H-uniqueness}) holds true, then the reflecting process $K$ is also unique. 
	Moreover, almost every path of $K$ is absolutely continuous with respect to the Lebesgue measure. 
\end{thm}

\section{Particle system and mean-field limit}\label{sec-limit}

\subsection{Definition and well-posedness of the particle system}

The objective of this section is to establish the propagation of chaos result for the mean-field reflected BSDE (\ref{MF-RBSDE})-(\ref{refl}). Let $\{B^i\}_{1\leq i\leq N}$ be $N$ independent $d$-dimensional standard Brownian motions and denote by $\mathbf{F}^i=(\Fc^i_t)_{0\leq t\leq T}$ the filtration generated by $B^i$. Let $\xi^i$ be independent $\Fc^i_T$-measurable copies of $\xi$. We will consider the following interacting particle system: for all $i\in\{1,2,\cdots,N\}$ and $t\in[0,T]$,
\begin{equation}\label{RBSDE-N}
	\begin{dcases}
		 & \begin{aligned}
			   \widehat{Y}^i_t
			   = &
			   \widehat{\xi}^i+\int_t^T f(s,\widehat{Y}^i_s,\widehat{Z}^{i,i}_s,\widehat{\nu}^{N}_{s})ds - \int_t^T \sum_{j=1}^N \widehat{Z}^{i,j}_s dB^j_s + \int_t^T \partial_y H(\widehat{Y}^{i}_{s},\widehat{\mu}^{N}_{s}) d\widehat{K}^i_s \\
			     & +\frac{1}{N}\sum_{j=1}^N \int_t^T \partial_\mu H(\widehat{Y}^{j}_{s},\widehat{\mu}^{N}_{s})(\widehat{Y}^{i}_{s}) d\widehat{K}^j_s ,
		   \end{aligned} \\
		 & H(\widehat{Y}^i_t, \widehat{\mu}^{N}_t)\geq 0,\quad \int_0^T H(\widehat{Y}^i_s, \widehat{\mu}^{N}_s)d\widehat{K}^i_s=0 ,
		\\
	\end{dcases}
\end{equation}
where $\widehat{\mu}^{N}_t:=N^{-1}\sum_{j=1}^N \delta_{\widehat{Y}^j_t}$ and $\widehat{\nu}^{N}_t:=N^{-1}\sum_{j=1}^N \delta_{(\widehat{Y}^j_t,\widehat{Z}^{j,j}_t)}$. 
For each $i,j\in\{1,2,\cdots,N\}$, $(\widehat{Y}^i,\widehat{Z}^{i,j},\widehat{K}^i)$ are progressively measurable processes, and $\widehat{K}^i$ is continuous and non-decreasing with $\widehat{K}^i_0=0$. 
Note that we have to define the terminal condition $\widehat{\xi}^i$ in the above particle system, so that $\widehat{\xi}^i$ is $\Fc^i_T$-measurable and satisfies the constraint $H(\widehat{\xi}^i,\widehat{\mu}^{N}_T)\geq0$. The following lemma constructs such modified terminal condition $\widehat{\xi}^i$ using $\xi^i$.

\begin{lem}\label{lem-eta-i}
	Let Assumption \ref{H} be satisfied. Then, there exist $C>0$ and $N$ independent $\widehat{\xi}^i$ is $\Fc^i_T$-measurable random variables $(\widehat{\xi}^i)_{1\leq i\leq N}$ with $\e|\widehat{\xi}^i|^4<\infty$, such that
	\begin{equation}\label{est-xi}
		\begin{aligned}
			H(\widehat{\xi}^i,\widehat{\mu}^{N}_T)\geq0
			\,\,\,\mathrm{and}\,\,\,
			\e\left| \widehat{\xi}^i-\xi^i \right|^2 \leq C \e \left[ W^2_2(\mu^N_T,\mu_T)\right], \quad i=1,2,\cdots,N, 
		\end{aligned}
	\end{equation}
	where $\widehat{\mu}^{N}_T=N^{-1}\sum_{j=1}^N \delta_{\widehat{\xi}^i}$, $\mu^N_T=N^{-1}\sum_{j=1}^N \delta_{\xi^i} $ and $\mu_T=\Lc(\xi)$.
\end{lem}

\begin{proof}
	We define $(X^i)_{1\leq i\leq N}$ the solution of
	\begin{equation}\label{e31}
		\begin{aligned}
			X^i_t=\xi^i+\int_0^t \partial_y H(X^{i}_{s},\mu^{X,N}_{s}) ds ,
		\end{aligned}
	\end{equation}
	where $\mu^{X,N}:=N^{-1}\sum_{j=1}^N \delta_{X^j}$. Then, from \eqref{H-beta}, \eqref{Hx-Hy} and \eqref{H1-H2}, we obtain 
	\begin{equation}\label{e17}
		\begin{aligned}
			& H(X^i_t,\mu^{X,N}_t) - H(\xi^i,\widehat{\mu}^{N}_T) \\
			=    
			& \int_0^t |\partial_y H(X^{i}_{s},\widehat{\mu}^{X,N}_{s}) |^2ds
			  +\frac{1}{N}\sum_{j=1}^{N} \int_0^t \partial_\mu H(X^{i}_{s},\mu^{X,N}_{s})(X^{j}_{s})\cdot  \partial_y H(X^j_s,\mu^{X,N}_s)ds 
			\\
			\geq &
			\beta^2 t.
			\\
		\end{aligned}
	\end{equation}
	Set $t^*:=\inf \{ t\geq0 : H(X^i_t,\mu^{X,N}_t) \geq0 \}$ and $\widehat{\xi}^i:=X^i_{t^*}$. Then, we have $t^*\leq \beta^{-2} H^-(\xi^i,\mu^N_T)$. From $H^-(\xi^i,\mu_T)=0$ and the Lipschitz continuity of $H$, we obtain
	\begin{equation}
		\begin{aligned}
			& \e\left| \widehat{\xi}^i-\xi^i \right|^2
			= \e\left| X^i_{t^*} - X^i_0 \right|^2 
			\leq M^2 \e\left(t^*\right)^2 
			\\
			\leq & 
			\frac{M^2}{\beta^4} \e \left| H^-(\xi^i,\mu^N_T) - H^-(\xi^i,\mu_T)  \right|^2
			\\
			\leq & 
			\frac{M^4}{\beta^4} \e \left[ W^2_2(\mu^N_T,\mu_T)\right].
		\end{aligned}
	\end{equation}
	The result follows.
\end{proof}

\begin{pr}
	Let Assumption \ref{H} be satisfied and (\ref{H-uniqueness}) holds true. 
	Then, the particle system (\ref{RBSDE-N}) is well posed. 
	Moreover, there exists $C>0$ independent of $N$, such that for all $1\leq i\leq N$, 
	\begin{equation}\label{est-5.2}
	\begin{gathered}
		\e\left[\sup_{0\leq t \leq T}	\left| Y^i_t \right|^2 + \int_0^T \sum_{j=1}^N\left|Z^{i,j}_s\right|^2ds +\int_0^T |k^i_s|^2 ds \right]\leq C , \\
		\sup_{0\leq t \leq T} \e \left| Y^i_t \right|^4 \leq C ,
	\end{gathered}
	\end{equation}
	where $k^i_t:=\frac{dK^i_t}{dt}$ is the Radon–Nikodym derivative of $K^i$ with respect to the Lebesgue measure. 
\end{pr}

\begin{proof}
	The uniqueness is obtained in an identical way to the proof of Proposition \ref{pr-unique}. 
	
	For $m\geq1$, we set 
	\begin{equation}
	\begin{aligned}
		\widehat{k}^{m,i}_t:=m H^-(\widehat{Y}^{m,i}_t,[\widehat{Y}^{m,i}_t]), \quad 
		 \widehat{K}^{m,i}_t:=\int_0^t \widehat{k}^{m,i}_s \,ds,
	\end{aligned}
	\end{equation}
	and define a penalized system of BSDE similar to \eqref{penal}. 
	Then, the existence is obtained by repeating the procedure in Section \ref{sec-exist}. 
\end{proof}

\subsection{The mean-field convergence}

In this section, we study the mean-field convergence of the solution of the particle system \eqref{RBSDE-N}. 

Let $(Y^i,Z^i,K^i)_{1\leq i\leq N}$ be the solution of the mean-field reflected BSDE (\ref{MF-RBSDE})-(\ref{refl}) with data $(B^i,\xi^i)_{1\leq i\leq N}$. Denote by $\Delta Y^i:=\widehat{Y}^i-Y^i$ and $\Delta Z^{i,j}:=\widehat{Z}^{i,j}-Z^i \delta_{i,j}$. 
Then, the following propagation of chaos result holds.

\begin{thm}\label{thm-limit}
	Let Assumption \ref{H} be satisfied and (\ref{H-uniqueness}) holds true. Then, there exists $C>0$, such that 
	\begin{equation}\label{est-poc}
		\begin{aligned}
			\e\left[\sup_{0\leq t \leq T} \frac{1}{N}\sum_{i=1}^{N}|\Delta Y^i_t|^2+\int_0^T \frac{1}{N}\sum_{i,j=1}^N |\Delta Z^{i,j}_s|^2ds
				\right]
			\leq
			C\left( \frac{1}{\sqrt{N}} +\sup_{0 \leq t \leq T} \e^{\frac{1}{8} }\left[W^4_2(\mu^N_t,\mu_t)\right]  \right).
		\end{aligned}
		\end{equation}
	\end{thm}
	
	\begin{proof}
		Let $\alpha$ be a large enough constant. For $1\leq i\leq N$, we apply \ito to $e^{\alpha t} |\Delta Y^i_t|^2$. From the Lipschitz continuity of $f$, we obtain 
		\begin{equation}\label{t5.3-e1}
		\begin{aligned}
			& e^{\alpha t} |\Delta Y^i_t|^2+\frac{1}{2} \int_t^T e^{\alpha s} \sum_{j=1}^{N} |\Delta Z^{i,j}_s|^2ds
			\\
		\leq &
			e^{\alpha T} |\Delta \xi^i|^2 		
			+ \int_t^T e^{\alpha s}\left( \frac{1}{N}\sum_{j=1}^{N} |\Delta Y^j_s|^2 + W^2_2(\mu^N_s,\mu_s) \right) ds
			-\int_t^T 2 e^{\alpha s} \Delta Y^i_s \cdot \sum_{j=1}^{N} \Delta Z^{i,j}_s dB^j_s \\
			& + J_1^i(t,T) + J_2^i(t,T),
		\end{aligned}
		\end{equation}
		with 
		\begin{equation}\label{eq:J1}
		\begin{aligned}
			J_1^i(t,T) & = \int_t^T 2 e^{\alpha s} \Delta Y^i_s \cdot \left( \partial_y H(\widehat{Y}^{i}_{s},\widehat{\mu}^{N}_{s}) \widehat{k}^i_s +\frac{1}{N}\sum_{j=1}^{N} \partial_\mu H(\widehat{Y}^{j}_{s},\widehat{\mu}^{N}_{s})(\widehat{Y}^{i}_{s}) \widehat{k}^{j}_{s}  \right) ds
		\end{aligned}
		\end{equation}
		and 
		\begin{equation}\label{eq:J2}
		\begin{aligned}
			J_2^i(t,T)&= -\int_t^T 2 e^{\alpha s} \Delta Y^i_s \cdot \left( \partial_y H(Y^{i}_{s},\mu_{s}) k^i_s + \widetilde{\e} \left[ \partial_\mu H(\widetilde{Y}^{}_{s},\mu^{}_{s})(Y^{i}_{s}) \widetilde{k}_s \right]  \right) ds . 
		\end{aligned}
		\end{equation}
		
	We sum the preceding inequality over $i$, take the expectation, and then examine each term on the r.h.s. separately. 

	1. We now deal with $J_1^i(t,T)$. From the concavity of $H$, we obtain 
	\begin{equation}\label{t5.3-ej1}
		\begin{aligned}
			& \frac{1}{N} \sum_{i=1}^{N} \e \left[ J_1^i(t,T) \right] 
			\\
		= & 
			\frac{1}{N}\sum_{i=1}^{N} \e \left[ \int_t^T 2 e^{\alpha s} \widehat{k}^i_s \left( 
			\Delta Y^i_s \cdot\partial_y H(\widehat{Y}^{i}_{s},\widehat{\mu}^{N}_{s}) 
			+ \frac{1}{N}\sum_{j=1}^{N} \Delta Y^j_s \cdot\partial_\mu H(\widehat{Y}^{i}_{s},\widehat{\mu}^{N}_{s})(\widehat{Y}^{j}_{s}) 
			\right) ds \right]
			\\
			& + \frac{1}{N^2}\sum_{i,j=1}^{N}\e \left[ 
				\int_t^T  2 e^{\alpha s} \left(- \widehat{k}^i_s \Delta Y^j_s \cdot\partial_\mu H(\widehat{Y}^{i}_{s},\widehat{\mu}^{N}_{s})(\widehat{Y}^{j}_{s})  
				+ \widehat{k}^j_s \Delta Y^i_s\cdot \partial_\mu H(\widehat{Y}^{j}_{s},\widehat{\mu}^{N}_{s})(\widehat{Y}^{i}_{s}) 
				\right) ds \right]
			\\
		\leq &
			\frac{1}{N}\sum_{i=1}^{N} \e \left[ \int_t^T e^{\alpha s} \left( H(\widehat{Y}^i_s,\widehat{\mu}^N_s)-H(Y^i_s,\mu^N_s) \right) \widehat{k}^i_s  ds \right]
			\\
		\leq & 
			\frac{1}{N}\sum_{i=1}^{N} \e \left[ \int_t^T e^{\alpha s} \left( H(Y^i_s,\mu_s)-H(Y^i_s,\mu^N_s) \right) \widehat{k}^i_s  ds \right]
			\\
		\leq &
			\frac{C}{N}\sum_{i=1}^{N} \e^{\frac{1}{2}} \left[ \int_t^T W_2^2(\mu^N_s,\mu_s) ds \right] 
			\,
			\e^{\frac{1}{2}} \left[ \int_t^T (\widehat{k}^i_s)^2 ds \right] 
		\\
		\leq &
			C \sup_{t \leq s \leq T} \e^{\frac{1}{2}}\left[W_2^2(\mu^N_s,\mu_s)\right]. 
		\\
		\end{aligned}
	\end{equation}

	2. The mean of the second term $J_2^i(t,T)$ reads 
	\begin{equation}\label{t5.3-l123}
		\begin{aligned}
			& \frac{1}{N}\sum_{i=1}^{N} \e \left[ J_2^i(t,T) \right] = I_1(t,T) + I_2(t,T) + I_3(t,T),  
			\\
		\end{aligned}
	\end{equation}
	where 
	\begin{equation}\label{eq:I1}
	\begin{aligned}
		I_1(t,T)&= \frac{2}{N}\sum_{i=1}^{N} \e \left[ \int_t^T  e^{\alpha s} k^i_s  \Delta Y^i_s \cdot \left(- \partial_y H(Y^{i}_{s},\mu^{}_{s})  +  \partial_y H(Y^{i}_{s},\mu^{N}_{s})  \right) ds \right] , 
	\end{aligned}
	\end{equation}
	\begin{equation}\label{eq:I2}
	\begin{aligned}
		I_2(t,T)&= -\frac{2}{N^2}\sum_{i,j=1}^{N} \e \left[ \int_t^T 
			 e^{\alpha s} k^i_s \left( \Delta Y^i_s \cdot \partial_y H(Y^{i}_{s},\mu^{N}_{s}) 
			+ \Delta Y^j_s \cdot \partial_\mu H(Y^{i}_{s},\mu^{N}_{s})(Y^{j}_{s}) 
			\right) 
			ds \right] , 
	\end{aligned}
	\end{equation}
	and 
	\begin{equation}\label{eq:I3}
	\begin{aligned}
		I_3(t,T)&= \frac{2}{N}\sum_{i=1}^{N} \e \left[ \int_t^T e^{\alpha s} \Delta Y^i_s \cdot\left(
				\partial_\mu H(Y^{j}_{s},\mu^{N}_{s})(Y^{i}_{s}) k^j_s
				- \widetilde{\e} \left[ \partial_\mu H(\widetilde{Y}^{}_{s},\mu^{}_{s})(Y^{i}_{s})\widetilde{k}_s  \right]
			\right) ds \right].
	\end{aligned}
	\end{equation}

	2.1. For the first part, from \eqref{est-5.2}, the Lipschitz continuity of $\partial_y H$ and Hölder's inequality, we obtain 
	\begin{equation}\label{t5.3-eqi1}
	\begin{aligned}
		I_1(t,T)
	\leq &
		\frac{C}{N}\sum_{i=1}^{N} \e \left[ \int_t^T  e^{\alpha s} k^i_s  |\Delta Y^i_s| W_2(\mu^N_s,\mu_s) ds \right]
	\\
	\leq & 
		\frac{C}{N}\sum_{i=1}^{N} \e^{\frac{1}{2} } \left[\int_t^T (k^i_s)^2 \,ds\right]
		\e^{\frac{1}{4} } \left[ \int_t^T |\Delta Y^i_s|^4 \,ds\right]
		\e^{\frac{1}{4} } \left[ \int_t^T W^4_2(\mu^N_s,\mu_s) \,ds\right]
	\\
	\leq &
		C \sup_{t \leq s \leq T} \e^{\frac{1}{4} } \left[ W^4_2(\mu^N_s,\mu_s) \right]. 
	\end{aligned}
	\end{equation}

	2.2. For the second part, from the concavity of $H$, the Skorokhod condition, and the inequality $H(\widehat{Y}^i_s,\mu^N_s)\geq0$, we obtain
	\begin{equation}\label{t5.3-eq92}
	\begin{aligned}
		I_2(t,T)
	\leq &
		\frac{2}{N}\sum_{i=1}^{N} \e \left[ \int_t^T \left(
			H(Y^i_s,\mu^N_s)-H(\widehat{Y}^i_s,\widehat{\mu}^N_s)
		\right) k^i_s ds \right]
	\\
	\leq &
		\frac{2}{N}\sum_{i=1}^{N} \e \left[ \int_t^T \left(
			H(\widehat{Y}^i_s,\mu^N_s)-H(\widehat{Y}^i_s,\widehat{\mu}^N_s)
		\right) k^i_s ds \right]
	\\
	\leq &
		C \sup_{t \leq s \leq T} \e^{\frac{1}{2}} \left[ W^2_2(\mu^N_s,\mu_s) \right]. 
	\\
	\end{aligned}
	\end{equation}
		
	2.3. Then, we deal with $I_3(t,T)$. We split it into two parts: 
	\begin{equation}
		\begin{aligned}
			I_3(t,T) = & (a) + (b),  \\
		\end{aligned}
		\end{equation}
		where 
		\begin{equation}\label{eq:(a)}
		\begin{aligned}
			(a) = \frac{2}{N^2} \sum_{i, j=1}^N \e\left[\int_t^T e^{\alpha s}\Delta Y_s^i \cdot \left( \partial_\mu H(Y_s^j, \mu_s^N)(Y_s^i) 
			-\partial_\mu H(Y_s^j, \mu_s)(Y_s^i) 
			\right)
			k_s^j d s\right]
		\end{aligned}
		\end{equation}
		and 
		\begin{equation}\label{eq:(b)}
		\begin{aligned}
			(b) = \frac{2}{N^2} \sum_{i, j=1}^N \e\left[\int_t^T e^{\alpha s} \Delta Y_s^i \cdot\left(\partial_\mu H(Y_s^j, \mu_s)(Y_s^i) k_s^j-\widetilde{\e}\left[\partial_\mu H(\widetilde{Y}_s, \mu_s)\left(Y_s^i\right) \widetilde{k}_s\right]\right) d s\right]. 
		\end{aligned}
		\end{equation}
		From the Lipschitz continuity of $\partial_{\mu} H$ and \eqref{est-5.2}, we obtain in an identical way to the calculation of \eqref{t5.3-eqi1}: 
		\begin{equation}\label{t5.3-eq94}
			\begin{aligned}
			(a) \leq C \sup_{t \leq s \leq T} \e^{\frac{1}{4} } \left[ W^4_2(\mu^N_s,\mu_s) \right]. \\
			\end{aligned}
		\end{equation}
		Recalling that $(Y^i)_{1\leq i\leq N}$ are independent of each other, we obtain 
		\begin{equation}
		\begin{aligned}
			\e\left[\int_t^T \Delta Y_s^i \cdot\left(\partial_\mu H(Y_s^j, \mu_s)(Y_s^i) k_s^j-\widetilde{\e}\left[\partial_\mu H(\widetilde{Y}_s, \mu_s)\left(Y_s^i\right) \widetilde{k}_s\right]\right)ds\right]=0,\quad \forall 1\leq i\neq j\leq N.  \\
		\end{aligned}
		\end{equation}
		From the preceding equation and the boundedness of $\partial_\mu H$, we obtain  
		\begin{equation}
		\begin{aligned}
		(b) =& 
		\frac{2}{N^2} \sum_{i=1}^N \e\left[\int_t^T e^{\alpha s} \Delta Y_s^i \cdot\left(\partial_\mu H(Y_s^i, \mu_s)(Y_s^i) k_s^i-\widetilde{\e}\left[\partial_\mu H(\widetilde{Y}_s, \mu_s)\left(Y_s^i\right) \widetilde{k}_s\right]\right) d s\right] 
		\\
		\leq &\frac{C}{N^2} \sum_{i=1}^N \e^{\frac{1}{2}}\left[\int_t^T \left|\Delta Y_s^i\right|^2 d s\right] \e^{\frac{1}{2}}\left[\int_t^T\left(\left(k_s^i\right)^2+\left(k_s\right)^2\right) d s\right]
		\leq \frac{C}{N}.\\
		\end{aligned}
		\end{equation}
	
		Bringing together the above estimates, we have 
		\begin{equation}\label{t5.3-ej23}
		\begin{aligned}
			&\frac{1}{N} \sum_{i=1}^N \e\left[J_2^i(t, T)\right] \leq C\left(\frac{1}{N}+\sup _{t \leq s \leq T} \e^{\frac{1}{4} }\left[W_2^4\left(\mu_s^N, \mu_s\right)\right]\right). 
		\end{aligned}
		\end{equation}
	
		From \eqref{est-5.2}, \eqref{t5.3-e1}, \eqref{t5.3-ej1}, and Gronwall's lemma, we obtain  
		\begin{equation}\label{t5.3-supe}
		\begin{aligned}
			&\sup _{0 \leq t \leq T} \e\left[\frac{1}{N} \sum_{i=1}^N\left|\Delta Y_t^i\right|^2\right]+\e\left[\int_0^T \frac{1}{N} \sum_{i, j=1}^N\left|\Delta Z_s^{i, j}\right|^2 d s\right]\\
			\leq& C\left(\frac{1}{N} \sum_{i=1}^N \e\left|\Delta \xi^i\right|^2
			+\sup _{0 \leq t \leq T} \e^{\frac{1}{2}}\left[W_2^2\left(\mu_t^N, \mu_t \right)\right]
			+\frac{1}{N}
			+\sup _{0 \leq t \leq T} \e^{\frac{1}{4} }\left[W_2^4\left(\mu_t^N, \mu_t\right)\right]\right)\\
			\leq& C\left(\frac{1}{N}+\sup _{0 \leq t \leq T} \e^{\frac{1}{4}}\left[W_2^4\left(\mu_t^N, \mu_t\right)\right]\right). \\
		\end{aligned}
		\end{equation}

	Coming back to \eqref{t5.3-e1}, summing the inequality over $i$, taking the supremum in $t$ and the expectation, we deduce from BDG inequality  
	\begin{equation}\label{t5.3-esup}
		\begin{aligned}
			&\e\sup _{0 \leq t \leq T}\left[ \frac{1}{N} \sum_{i=1}^N e^{\alpha t}\left|\Delta Y_t^i\right|^2
			\right]
			+\e\left[\int_0^T \frac{1}{N} \sum_{i, j=1}^N e^{\alpha s}\left|\Delta Z_s^{i, j}\right|^2 d s\right]\\
			\leq &\frac{C}{N} \e \left[ \sum_{i=1}^N e^{\alpha T}\left|\Delta\xi^i\right|^2
			+ \int_t^T e^{\alpha s} |\Delta Y^i_s|^2  \,ds
			+ \sum_{k=1}^2 \sup _{0 \leq t \leq T}\sum_{i=1}^N J_k^i(t, T)
			\right]. \\
		\end{aligned}
	\end{equation}

	In view of the definition \eqref{t5.3-e1} of $J^i_k(t,T)$, from the boundedness of $\partial_y H$ and $\partial_\mu H$, we obtain 
	\begin{equation}\label{t5.3-esup-j}
	\begin{aligned}
		& \frac{1}{N} \e\left[\sum_{k=1}^2\sup _{0 \leq t \leq T}\sum_{i=1}^N J_k^i(t, T)\right] \\
	\leq &
		\frac{C}{N} \sum_{i=1}^N \e\left[\int_0^T e^{\alpha s} \left(k_s^i+\e[k_s^i]\right) \left|\Delta Y_s^i\right| d s\right]\\
	\leq &
		\frac{C}{N} \sum_{i=1}^N \e^{\frac{1}{2} }\left[\int_0^T \left(k_s^i\right)^2 d s\right]
		\e^{\frac{1}{2} }\left[\int_0^T \left|\Delta Y_s^i\right| ^2 ds\right] \\
	\leq &
		C\left(\frac{1}{\sqrt{N}}+\sup _{0 \leq t \leq T} \e^{\frac{1}{8}}\left[W_2^4\left(\mu_t^N, \mu_t\right)\right]\right).\\
	\end{aligned}
	\end{equation}
	From \eqref{est-xi}, \eqref{t5.3-esup}, the last two inequalities, and Gronwall's lemma, we obtain the desired result. 
\end{proof}

\begin{cor}
	Let Assumption \ref{H} be satisfied and (\ref{H-uniqueness}) holds true. Assume that $f^0\in H^{q,n}$ and $\e|\xi|^q<\infty$ for some $q>4$. Then, there exists $C>0$, such that 
	\begin{equation}\label{est-poc-cor}
	\begin{aligned}
		\e\left[\sup_{0\leq t \leq T} \frac{1}{N}\sum_{i=1}^{N}|\Delta Y^i_t|^2+\int_0^T \frac{1}{N}\sum_{i,j=1}^N |\Delta Z^{i,j}_s|^2ds
			\right]
		\leq C 
		\begin{dcases}
			N^{-\frac{1}{8} }, & n<4, \\
			N^{-\frac{1}{8} } \ln (N+1) , & n=4, \\
			N^{-\frac{1}{2n} } , & n>4. \\
		\end{dcases}
	\end{aligned}
	\end{equation}
\end{cor}

\begin{proof}
	By repeating the proof of Lemma \ref{lem-a-priori-supE}, we can show that 
	\begin{equation}
	\begin{aligned}
		\sup_{0 \leq t \leq T} \e|\widehat{Y}^i_t|^q<\infty, \quad i=1,2,\cdots,N. 
	\end{aligned}
	\end{equation}
	Hence, we obtain the result by combining the preceding theorem and \cite[Theorem 1]{FG15}. 
\end{proof}

\section{Related obstacle problem for PDEs in Wasserstein space}\label{sec:obstacle problem}

In this section, we connect the MF-RBSDE with an obstacle problem for PDEs in Wasserstein space. Throughout this section, we assume that $Y_t$ is a one-dimensional process, i.e., $n=1$. Note that \eqref{H-beta} and \eqref{Hx-Hy} imply that $\partial_y H(y^{}_{},\mu^{}_{})$ and $\partial_\mu H(y^{}_{},\mu^{}_{})$ are either always positive or always negative in the one-dimensional case. Without loss of generality, we assume that 
\begin{equation}\label{eq:Hy,Hmu}
\begin{aligned}
    \partial_y H(y^{}_{},\mu^{}_{}) < 0 \,\,\, \mathrm{and} \,\,\, \partial_\mu H(y^{}_{},\mu^{}_{})(v)<0, \quad \forall (y,v,\mu)\in \R\times\R\times\Pc_2(\R). 
\end{aligned}
\end{equation}

We consider a Markovian setup, where the terminal condition $\xi$ is given by $g(X^{t, X_0}_T, [X^{t, X_0}_T])$, and $X^{t, X_0}$ is the solution of an SDE. Namely, we consider the following Forward-Backward SDE with mean-field reflection (MF-RFBSDE): 
\begin{equation}\label{eq:FBSDE}
\begin{dcases}
\begin{aligned}
    X^{t, X_0}_s =& X_0 + \int_t^s b(r, X^{t, X_0}_r) \,dr + \int_t^s \sigma(r, X^{t, X_0}_r) \,dB_r, \quad s\in[t,T], \\
    X^{t, X_0}_s =& X_0, \quad s\in[0,t)  \\
    Y^{t, X_0}_s =& g(X^{t, X_0}_T, [X^{t, X_0}_T]) + \int_s^T f(r, X^{t, X_0}_r, Y^{t, X_0}_r, [X^{t, X_0}_r], [Y^{t, X_0}_r]) \,dr \\
    &- \int_s^T Z^{t, X_0}_r \,dB_r + R^{t, X_0}_T - R^{t, X_0}_s , \quad s\in[t,T], \\
    Y^{t, X_0}_s =& Y^{t, X_0}_t, \quad s\in[0,t), \\
        H( Y^{t,  X_0 }_s & ,  [Y^{t, X_0}_s])  \geq0, \quad s\in[t,T], \quad
    \int_t^T H(Y^{t, X_0}_r,[Y^{t, X_0}_r]) \,dK^{t, X_0}_r = 0, \\
\end{aligned}
\end{dcases}
\end{equation}
where 
\begin{equation}\label{eq:R-def}
\begin{aligned}
    & R^{t, X_0}_s := \int_t^s \partial_y H(Y^{t, X_0}_{r}, [Y^{t, X_0}_r]) \,dK^{t, X_0}_r 
    + \widetilde{\e} \left[ \int_t^s  \partial_\mu H(\widetilde{Y}^{t, X_0}_{r}, [Y^{t, X_0}_r])(Y^{t, X_0}_{r}) \, d \widetilde{K}^{t, X_0}_r \right], & s\in[t,T], 
\end{aligned}
\end{equation}
and the superscript $(t,X_0)$ stands for the initial condition of the SDE. 

\begin{rem}\label{rem:f-z-mu}
	We assume that $f$ does not depend on $Z$, and only depends on the marginal distribution $\mu$ of $\nu$. This assumption is due to the lack of a comparison principle for mean-field FBSDEs (see \cite[Section 3]{BLP09}), which poses difficulty in analyzing the existence of a decoupling field when the driver depends on the $Z$ argument.
\end{rem}

The coefficients $b, \sigma$ of the SDE and the functions $f,g$ satisfy the following conditions.
\begin{Assumption}\label{H:SDE}
\begin{enumerate}
    \item\label{H:b,sigma} The functions $b$ and $\sigma$ are progressively measurable mappings from $\Omega \times [0, T] \times \R^l$ to $\R^l$ and $\R^{l \times d}$, respectively, such that for all $s\in[0,T]$ and $x_1, x_2 \in \R^l$,
    \begin{equation}\label{eq:b,sigma}
    \begin{gathered}
        \left| b(s, x_1) - b(s, x_2) \right| +\left| \sigma(s, x_1) - \sigma(s, x_2) \right| \leq L \left| x_1 - x_2 \right|, \\
        \left| b(s, 0) \right| + \left| \sigma(s, 0) \right| \leq L . 
    \end{gathered}
    \end{equation}
    \item The driver $f$ is a mapping from $[0,T]\times \R^l \times \R \times \Pc_2(\R^l) \times \Pc_2(\R)$ to $\R$, which satisfies the linear growth condition and the Lipschitz condition, i.e., there exists $L>0$, such that for all $t\in[0,T], x,x_1,x_2\in\R^l, y_1,y_2\in\R, \lambda,\lambda_1,\lambda_2 \in \Pc_2(\R^l)$ and $\mu_1,\mu_2 \in \Pc_2(\R)$, 
    \begin{equation}\label{eq:f-linear-growth}
    \begin{aligned}
        \left| f(t,x,0,\lambda,\delta_0) \right| \leq L \left( 1+ \left| x \right| + W_2(\lambda,\delta_0) 
		 \right),  
    \end{aligned}
    \end{equation}
    \begin{equation}\label{eq:f-Lipschitz-2}
    \begin{aligned}
       & \left| f(t,x_1,y_1,\lambda_1,\mu_1) - f(t,x_2,y_2,\lambda_2,\mu_2) \right|  \\
    \leq &
        L\left(|x_1-x_2| + |y_1-y_2| +W_2(\lambda_1,\lambda_2) +W_2(\mu_1,\mu_2) \right).
    \end{aligned}
    \end{equation}
	\item\label{H:g} The function $g$ is a Lipschitz continuous mapping from $\R^l \times \Pc_2(\R)$ to $\R$, i.e., for all $x_1,x_2\in\R^l$ and $\lambda_1,\lambda_2 \in \Pc_2(\R^l)$, 
	\begin{equation}\label{eq:g-Lipschitz}
	\begin{aligned}
		& \left| g(x_1,\lambda_1) - g(x_2,\lambda_2) \right|  
		\leq 
			L\left(|x_1-x_2| +W_2(\lambda_1,\lambda_2) \right). 
	\end{aligned}
	\end{equation} 
\end{enumerate}
\end{Assumption}

Under these assumptions, the SDE part of \eqref{eq:FBSDE} has a unique strong solution $X^{t, X_0}$. From Theorem \ref{thm-exist-unique}, there exists a unique solution $(Y^{t, X_0},Z^{t, X_0},K^{t, X_0})$ to the mean-field reflected BSDE part of \eqref{eq:FBSDE}. Therefore, the MF-RFBSDE \eqref{eq:FBSDE} admits a unique solution. 

We will connect the solution of \eqref{eq:FBSDE} to a viscosity solution of an obstacle problem for PDE. Consider the following problem in Wasserstein space:
\begin{equation}\label{PDE}
\begin{aligned}
\begin{dcases}
    \min \left\{ \left( \partial_{t} + \Lc_t \right) u(t,x,\lambda) + f(t,x,u(t,x,\lambda), \lambda, u(t,\cdot,\lambda)_* \lambda) , \,\,
    H(u(t,x,\lambda), u(t,\cdot,\lambda)_* \lambda) \right\} = 0, \\
    u(T,\cdot,\cdot) = g, \\
\end{dcases}
\end{aligned}
\end{equation}
where $u(t,\cdot,\lambda)_* \lambda$ stands for the pushforward measure of $\lambda$ under the mapping $u(t, \cdot, \lambda)$, and $\Lc_t$ is an operator given by 
\begin{equation}\label{eq:Lc}
\begin{aligned}
    \Lc_t \varphi(t,x,\lambda):=
    & b(t,x) \nabla_x \varphi(t,x,\lambda)  + \frac{1}{2} \mathrm{Tr}\left[ \left( \sigma\sigma^* \right)(t,x) \nabla^2_x \varphi(t,x,\lambda) \right] \\
    & + \e_{X\sim\lambda}\left[  
        b(t,X) \cdot \partial_\mu \varphi(t,x,\lambda)(X) 
        +  \frac{1}{2} \mathrm{Tr} \left[ \left( \sigma\sigma^* \right)(t,X) \partial_v \partial_\mu \varphi(t,x,\lambda)(X) \right]
     \right]
    , 
\end{aligned}
\end{equation}
for all smooth $\varphi:[0,T]\times \R^l \times \Pc_2(\R^l) \to \R $.

We define the notion of viscosity solution of \eqref{PDE} as follows: 
\begin{defn}\label{defn:viscosity}
    A function $u:[0,T]\times \R^l \times \Pc_2(\R^l) \to \R$ is called a viscosity subsolution (resp. supersolution) of the obstacle problem \eqref{PDE} if: 
    \begin{enumerate}
        \item the function $u$ is continuous and locally bounded, 
        \item for any $(t,x,\lambda)\in [0,T]\times \R^l \times \Pc_2(\R^l)$, for any test function $\varphi:[0,T]\times \R^l \times \Pc_2(\R^l) \to \R$ 
        (see \cite[Definition 11.18]{CD18} for the definition of test functions) 
        such that $u-\varphi$ has a global minimum (resp. maximum) in $(t,x,\lambda)$, we have 
        \begin{equation}\label{eq:sub-super solution}
        \begin{aligned}
            \min \left\{ \left( \partial_{t} + \Lc_t \right) \varphi(t,x,\lambda) + f(t,x,u(t,x,\lambda), u(t,\cdot,\lambda)_* \lambda) , \,\,
            H(u(t,x,\lambda), u(t,\cdot,\lambda)_* \lambda) \right\}
        \leq 0 \, (\mathrm{resp. } \geq0), 
        \end{aligned}
        \end{equation}
        \item $u(T,\cdot,\cdot)=g$ on $ \R^l \times \Pc_2(\R^l)$. 
    \end{enumerate}
    A function $u:[0,T]\times \R^l \times \Pc_2(\R^l) \to \R$ is called a viscosity solution if and only if it is a subsolution and a supersolution. 
\end{defn}

For any $(t,x,\lambda) \in [0,T]\times \R^l \times \Pc_2(\R^l)$, we define 
\begin{equation}\label{eq:decoupling}
\begin{aligned}
    u(t,x,\lambda) := Y^{t,x,\lambda}_t, 
\end{aligned}
\end{equation}
where $(X^{t,x,\lambda},Y^{t,x,\lambda},Z^{t,x,\lambda})$ is the solution of 
\begin{equation}\label{eq:FBSDE-x}
\begin{dcases}
\begin{aligned}
	X^{t, x, \lambda}_s =& x + \int_t^s b(r, X^{t, x, \lambda}_r) \,dr + \int_t^s \sigma(r, X^{t, x, \lambda}_r) \,dB_r, \quad s\in[t,T], \\
	X^{t, x, \lambda}_s =& x, \quad s\in[0,t). \\
	Y^{t, x, \lambda}_s =& g(X^{t, x, \lambda}_T, [X^{t, \lambda}_T]) + \int_s^T f(r, X^{t, x, \lambda}_r, Y^{t, x, \lambda}_r, [X^{t, \lambda}_r], [Y^{t,\lambda}_r]) \,dr \\
	&- \int_s^T Z^{t, x, \lambda}_r \,dB_r + R^{t, x, \lambda}_T - R^{t, x, \lambda}_s , \quad s\in[t,T], \\
	Y^{t, x, \lambda}_s =& Y^{t, x, \lambda}_t, \quad s\in[0,t), \\
\end{aligned}
\end{dcases}
\end{equation}
and 
\begin{equation}\label{eq:Rx-def}
\begin{aligned}
    &  R^{t, x, \lambda}_s := \e \left[ R^{t, \lambda}_s \mid X^{t, \lambda}_t = x \right], \quad s\in[t,T] . \\
\end{aligned}
\end{equation}

We present the following result that connects $u$ with the solution of the obstacle problem \eqref{PDE}. 
\begin{thm}\label{thm-viscosity}
    Under Assumption \ref{H} 2.(a)-(g) and \ref{H:SDE}, the decoupling field $u$ defined in \eqref{eq:decoupling} is a viscosity solution of the obstacle problem \eqref{PDE}. 
\end{thm}

\begin{rem}\label{rem:thm-viscosity}
    Theorem \ref{thm-viscosity} establishes the existence of a viscosity solution of the obstacle problem \eqref{PDE} in the Wasserstein space setting. However, this result does not address the issue of uniqueness for the viscosity solution. The study of uniqueness for viscosity solutions of PDEs in Wasserstein space is inherently difficult, as traditional techniques used for proving uniqueness, such as comparison principles, often do not apply in this context \cite[Section 3]{BLP09}.

    Recently, Talbi et al. \cite{TTZ22} introduced an alternative notion of viscosity solutions for obstacle problems in Wasserstein space, establishing the existence, uniqueness, stability, and comparison principles. While their results do not directly apply to our specific problem, their work may provide insights for future research. 
\end{rem}

The following lemma addresses the (joint) continuity of the decoupling field $u$. 

\begin{lem}\label{lem:u-continuous}
    Under Assumption \ref{H} 2.(a)-(g) and \ref{H:SDE}, the decoupling field $u$ defined in \eqref{eq:decoupling} is continuous in $[0,T]\times \R^l \times \Pc_2(\R^l)$.
\end{lem}
    
\begin{proof}
    For simplicity, we denote $X^i:=X^{t^i,x^i,\lambda^i}$ and $Y^i:=Y^{t^i,x^i,\lambda^i}$ for $i\geq0$. Fix $(t^0,x^0,\lambda^0)$ in $[0,T]\times \R^l \times \Pc_2(\R^l)$. For any sequence $(t^m,x^m,\lambda^m)_{m\geq1}$ converging to $(t^0,x^0,\lambda^0)$ as $m\to\infty$, we will prove that 
    \begin{equation}\label{eq:u-continuous}
    \begin{aligned}
    \left| u(t^m,x^m,\lambda^m) - u(t^0,x^0,\lambda^0) \right|
            = \left| Y^{m}_{t^m} - Y^{0}_{t^0} \right|
    \to 0.
    \end{aligned}
    \end{equation}
    
    We split it into two terms and then deal with them separately. Note that $Y^m_{t^m}$ and $Y^0_{t^0}$ are deterministic, so 
    \begin{equation}\label{eq:Ym-Y0}
    \begin{aligned}
        & \left| Y^m_{t^m} - Y^0_{t^0} \right| 
    =
    \left| \e \left[ Y^m_{t^m} - Y^0_{t^0} \right] \right| 
    \leq 
        \left| \e \left[ Y^m_{t^m} - Y^0_{t^m} \right] \right| 
        + \left| \e \left[ Y^0_{t^m} - Y^0_{t^0} \right] \right| .  \\
    \end{aligned}
    \end{equation}
    
    From Proposition \ref{pr:stability}, we have 
    \begin{equation}\label{eq:Ym-Y0:1}
    \begin{aligned}
        & \left| \e \left[ Y^m_{t^m} - Y^0_{t^m} \right] \right| 
    \leq 
        C \e^{\frac{1}{2} } \left[ \left| \Delta g \right|^2 + \int_0^T \left| \delta f(s,Y^0_s,[Y^{t^0,\lambda^0}_s]) \right|^2 \,ds \right], 
    \end{aligned}
    \end{equation}
    where 
    \begin{equation}\label{eq:delta-g-f}
    \begin{gathered}
        \Delta g := g(X^{m}_T, [X^{t^m, \lambda^m}_T]) - g(X^{0}_T, [X^{t^0, \lambda^0}_T]), \\
        \delta f(s,y,\mu):=f(s,X^m_s,y,[X^{t^m, \lambda^m}_s],\mu) - f(s,X^0_s,y,[X^{t^0, \lambda^0}_s],\mu).
    \end{gathered}
    \end{equation}
    From the Lipschitz continuity of $g$ and $f$ and the standard estimate \cite[Theorem 3.4.3]{Z17} of solutions of SDE, the r.h.s. of \eqref{eq:Ym-Y0:1} converges to $0$ as $m\to\infty$.
    
    Next, we will deal with the second term in \eqref{eq:Ym-Y0}. 
    \begin{equation}\label{eq:Ym-Y0:2}
    \begin{aligned}
        & \left| \e \left[ Y^0_{t^m} - Y^0_{t^0} \right] \right| 
    \leq 
         \e \left| \int_{t^0}^{t^m}   f(r,X^0_r,Y^0_r,[X^{t^0, \lambda^0}_r],[Y^{t^0,\lambda^0}_r])   \,dr \right| 
        + \e \left| R^0_{t^m} - R^0_{t^0} \right|. 
    \end{aligned}
    \end{equation}
    In view of the definition \eqref{eq:R-def} and \eqref{eq:Rx-def} of $R^0$, and \eqref{H-M}, we have 
    \begin{equation}\label{eq:E|R|}
	\begin{aligned}
		& \e \left| R^0_{t^m} - R^0_{t^0} \right| 
	\leq 
		C \e \left| \int_{t^0}^{t^m}  k^0_r \,dr \right| ,
	\end{aligned}
	\end{equation}
    where $k^0_r:=\frac{dK^0_r}{dr}$. Recall that we have $k^0\in L^2(\Omega\times[0,T])$ from the proof of Proposition \ref{pr-exist}. 

    From the growth condition of $f$ and the last two inequalities, we obtain 
    \begin{equation}\label{eq:Ym-Y0:2.1}
    \begin{aligned}
        & \left| \e \left[ Y^0_{t^m} - Y^0_{t^0} \right] \right| 
    \leq 
        C  
        \e \left| \int_{t^0}^{t^m}  \left( 1 + \left| X^0_r \right| + \left| Y^0_r \right| + \e \left| X^{t^0,\lambda^0}_r \right| +  \e \left| Y^{t^0,\lambda^0}_r \right| \right) \,dr \right|  
        + C \e \left| \int_{t^0}^{t^m}  k^0_r \,dr \right|  , 
    \end{aligned}
    \end{equation}
    which converges to $0$ as $m\to\infty$. 
    
\end{proof}

\begin{proof}[Proof of Theorem \ref{thm-viscosity}]\label{proof:thm-viscosity}
    Our approach is inspired by the proofs of \cite[Theorem 8.5]{EKPPM97} and \cite[Theorem 37]{BCDH20}. We will first prove that the function $u$ defined in \eqref{eq:decoupling} is a subsolution of the obstacle problem \eqref{PDE}. Let $(t,x,\lambda)\in [0,T]\times \R^l \times \Pc_2(\R^l)$. For any test function $\varphi$, such that $u-\varphi$ has a global minimum at $(t,x,\lambda)$, we aim to show that \eqref{eq:sub-super solution} is valid. 
    
    Recall that $u(t,x,\lambda) = Y^{t,x,\lambda}_t$. From the dynamics \eqref{eq:FBSDE} of $Y^{t,x,\lambda}$, we obtain, for all $s\in[t,T]$, 
    \begin{equation}\label{eq:Eu}
    \begin{aligned}
        & \e \left[ u(s,X^{t,x,\lambda}_s,[X^{t,\lambda}_s]) \right] \\
    = &
        u(t,x,\lambda)-\e \left[ \int_t^s f(r,X^{t,x,\lambda}_r,Y^{t,x,\lambda}_r,[X^{t,\lambda}_r],[Y^{t,\lambda}_r]) \,dr \right] - \e \left[ R^{t,x,\lambda}_s - R^{t,x,\lambda}_t \right].  \\
    \end{aligned}
    \end{equation}
    
    Applying Itô's formula, we derive
    \begin{equation}\label{eq:Ephi}
    \begin{aligned}
        \e \left[ \varphi(s,X^{t,x,\lambda}_s,[X^{t,\lambda}_s]) \right] 
    =
        \varphi(t,x,\lambda) + \e \left[ \int_t^s \left( \partial_{t} + \Lc_t \right) \varphi(r,X^{t,x,\lambda}_r,[X^{t,\lambda}_r])  \,dr \right] . 
    \end{aligned}
    \end{equation}
        
    Without loss of generality, we assume that $u(t,x,\lambda)=\varphi(t,x,\lambda)$. Since $u-\varphi$ is minimized at $(t,x,\lambda)$, we have 
    \begin{equation}\label{eq:u>=phi}
    \begin{aligned}
        \e \left[ u(s,X^{t,x,\lambda}_s,[X^{t,\lambda}_s]) \right] \geq \e \left[ \varphi(s,X^{t,x,\lambda}_s,[X^{t,\lambda}_s]) \right]. 
    \end{aligned}
    \end{equation}

    Combining \eqref{eq:Eu} with \eqref{eq:Ephi}, we deduce
    \begin{equation}\label{eq:E(u-phi)}
    \begin{aligned}
        & \e \left[ \int_t^s \left( \left( \partial_{t} + \Lc_t \right) \varphi(r,X^{t,x,\lambda}_r,[X^{t,\lambda}_r]) + f(r,X^{t,x,\lambda}_r,Y^{t,x,\lambda}_r,[X^{t,\lambda}_r],[Y^{t,\lambda}_r]) \right)  \,dr \right] \\
    \leq &
         - \e \left[ R^{t,x,\lambda}_s - R^{t,x,\lambda}_t \right] . 
    \end{aligned}
    \end{equation}

We now verify \eqref{eq:sub-super solution}. Assume now that $H(u(t,x,\lambda), u(t,\cdot,\lambda)_* \lambda))>0$. Then, from the continuity of $u$, there exists $\delta>0$, such that for all $(s,x',\lambda') \in [t,T]\times \R^l \times \Pc_2(\R^l)$ with $(s-t) + |x-x'| + W_2(\lambda,\lambda') \leq 3\delta$, we have 
\begin{equation}\label{eq:H>0}
\begin{aligned}
    H(u(s,x',\lambda'), u(s,\cdot,\lambda')_* \lambda')>0. 
\end{aligned}
\end{equation}
Moreover, from the standard estimate \cite[Theorem 3.4.3]{Z17} of solutions of SDE, there exists $s_1>t$, such that $W_2(\lambda,[X^{t,\lambda}_s])\leq \delta,\,\, \forall s\in [t,s_1]$. 

    Assume now that $s\in [t,s_1 \wedge (t+\delta)]$, we will estimate the r.h.s. of \eqref{eq:E(u-phi)}. 
	Then, we deduce from the Skorokhod condition that 
	\begin{equation}\label{eq:}
	\begin{aligned}
		\int_{t}^{s} \mathbf{1}_{\left\{ \left| X^{t,x,\lambda}_r - x \right| \leq \delta \right\}} k^{t,x,\lambda}_r \,dr = 0. 
	\end{aligned}
	\end{equation}
	From \eqref{eq:E|R|}, Young's inequality and Chebyshev's inequality, we deduce that 
    \begin{equation}\label{eq:ER}
    \begin{aligned}
        & \left| \e \left[ R^{t,x,\lambda}_s - R^{t,x,\lambda}_t \right] \right| 
    \leq 
        C \e  \int_{t}^{s}  k^{t,x,\lambda}_r \,dr 
    = 
        C \e  \int_{t}^{s} \mathbf{1}_{\left\{ \left| X^{t,x,\lambda}_r - x \right| > \delta \right\}} k^{t,x,\lambda}_r \,dr \\
    \leq &
        C \P^{\frac{1}{2} } \left[ \sup_{t \leq r \leq s}  \left| X^{t,x,\lambda}_r - x \right| > \delta  \right] 
        \e^{\frac{1}{2}} \left[ \int_{t}^{s} \left( k^{t,x,\lambda}_r \right)^2 \,dr \right] \\
    \leq & 
        \frac{C}{\delta^2} \e^{\frac{1}{2}} \left[  \sup_{t \leq r \leq s}  \left| X^{t,x,\lambda}_r - x \right|^4 \right]
        \e^{\frac{1}{2}} \left[ \int_{t}^{s} \left( k^{t,x,\lambda}_r \right)^2 \,dr \right] \\
    \leq &
        \frac{C}{\delta^2} (s-t) \e^{\frac{1}{2}} \left[ \int_{t}^{s} \left( k^{t,x,\lambda}_r \right)^2 \,dr \right] \\
    \end{aligned}
    \end{equation}

Dividing both sides of \eqref{eq:E(u-phi)} by $(s-t)$ and letting $s\to t$, from the last inequality and 
\begin{equation}\label{eq:Ek->0}
\begin{aligned}
    \lim_{s\to t} \e \left[ \int_{t}^{s} \left( k^{t,x,\lambda}_r \right)^2 \,dr \right] = 0, 
\end{aligned}
\end{equation}
we obtain 
\begin{equation}\label{eq:d_t phi}
\begin{aligned}
    \left( \partial_{t} + \Lc_t \right) \varphi(t,x,\lambda) + f(t,x,u(t,x,\lambda), u(t,\cdot,\lambda)_* \lambda) \leq 0. 
\end{aligned}
\end{equation}

Therefore, we have 
\begin{equation}\label{eq:subsolution}
\begin{aligned}
    \min \left\{ \left( \partial_{t} + \Lc_t \right) \varphi(t,x,\lambda) + f(t,x,u(t,x,\lambda), u(t,\cdot,\lambda)_* \lambda) , \,\,
        H(u(t,x,\lambda), u(t,\cdot,\lambda)_* \lambda) \right\}
    \leq 0. 
\end{aligned}
\end{equation}

We conclude the proof by showing that $u$ is also a supersolution of the obstacle problem \eqref{PDE}. We already know that $H(u(t,x,\lambda), u(t,\cdot,\lambda)_* \lambda) \geq 0$. Let $(t,x,\lambda) \in [0,T] \times \R^l \times \Pc_2(\R^l) \to \R$. For any test function $\varphi$, such that $u - \varphi$ has a global maximum at $(t,x,\lambda)$, we obtain the following using the same technique as in the calculation of \eqref{eq:E(u-phi)}: 
\begin{equation}\label{eq:E(u-phi),2}
\begin{aligned}
    & \e \left[ \int_t^s \left( \left( \partial_{t} + \Lc_t \right) \varphi(r,X^{t,x,\lambda}_r,[X^{t,\lambda}_r]) + f(r,X^{t,x,\lambda}_r,Y^{t,x,\lambda}_r,[X^{t,\lambda}_r],\mu^{t,\lambda}_r) \right)  \,dr \right] \\
    \geq &
    - \e \left[ R^{t,x,\lambda}_s - R^{t,x,\lambda}_t \right] 
    \geq 
    0,  \quad \forall s\in[t,T],  \\
\end{aligned}
\end{equation}
with the last inequality coming from \eqref{eq:Hy,Hmu} and the definitions \eqref{eq:R-def} and \eqref{eq:Rx-def} of $R^{t,x,\lambda}$. Dividing both sides of \eqref{eq:E(u-phi)} by $(s-t)$ and letting $s \to t$, we obtain
\begin{equation}\label{eq:d_t phi,2}
\begin{aligned}
\left( \partial_{t} + \Lc_t \right) \varphi(t,x,\lambda) + f(t,x,u(t,x,\lambda), u(t,\cdot,\lambda)_* \lambda) \geq 0.
\end{aligned}
\end{equation}

\end{proof}

\section*{Declarations}

This work was supported by National Natural Science Foundation of China (Grant No. 12031009). 
The author has no competing interests to disclose.

\end{document}